\documentclass[10pt]{article}

\usepackage{amsmath}
\usepackage{amsfonts}
\usepackage{amssymb}

\newtheorem{proposition}{Proposition}[section]
\newtheorem{theorem}{Theorem}[section]
\newtheorem{lemma}[theorem]{Lemma}
\newtheorem{corollary}[theorem]{Corollary}
\newtheorem{remark}[theorem]{Remark}

\newtheorem{definition}{Definition}

\def\phi{{\varphi}}

\DeclareSymbolFont{AMSb}{U}{msb}{m}{n}
\DeclareMathSymbol{\N}{\mathbin}{AMSb}{"4E}
\DeclareMathSymbol{\Z}{\mathbin}{AMSb}{"5A}
\DeclareMathSymbol{\R}{\mathbin}{AMSb}{"52}
\DeclareMathSymbol{\Q}{\mathbin}{AMSb}{"51}
\DeclareMathSymbol{\I}{\mathbin}{AMSb}{"49}
\DeclareMathSymbol{\C}{\mathbin}{AMSb}{"43}

\begin{document}

\title{Bessel potentials and optimal Hardy and Hardy-Rellich inequalities}
\author{ Nassif  Ghoussoub\thanks{Partially supported by a grant
from the Natural Sciences and Engineering Research Council of Canada.  } \quad  and \quad Amir %%@
Moradifam \thanks{Partially supported by a UBC Graduate Fellowship.  }
\\
\small Department of Mathematics,
\small University of British Columbia, \\
\small Vancouver BC Canada V6T 1Z2 \\
\small {\tt nassif@math.ubc.ca} \\
\small {\tt a.moradi@math.ubc.ca}
\\
%\today\\
%\date{January 20, 2005}\\
}
\maketitle

 \begin{abstract} We give necessary and sufficient conditions  on a pair of positive radial %%@
functions  $V$ and $W$ on a ball $B$ of radius $R$  in $R^{n}$, $n \geq 1$, so that the %%@
following inequalities hold for all $u \in C_{0}^{\infty}(B)$:   
\begin{equation*} \label{one}
\hbox{$\int_{B}V(x)|\nabla u |^{2}dx \geq \int_{B} W(x)u^2dx$,  }
\end{equation*} 
and
\begin{equation*} \label{two}
\hbox{$\int_{B}V(x)|\Delta u |^{2}dx \geq  \int_{B} W(x)|\nabla  %%@
u|^{2}dx+(n-1)\int_{B}(\frac{V(x)}{|x|^2}-\frac{V_r(|x|)}{|x|})|\nabla u|^2dx$.}
\end{equation*} 
This characterization makes a very useful connection between Hardy-type inequalities and the oscillatory behaviour of certain ordinary differential equations, and helps in the identification of  a large number of such couples $(V, W)$ -- that we call Bessel pairs -- as well %%@
as the best constants in the corresponding inequalities.  This allows us to   improve, extend, and %%@
unify many results --old and new-- about Hardy and Hardy-Rellich type inequalities,  such as those obtained by %%@
Caffarelli-Kohn-Nirenberg \cite{CKN}, Brezis-V\'{a}zquez \cite{BV},  Wang-Willem \cite{WW}, Adimurthi-Chaudhuri-Ramaswamy \cite{ACR}, Filippas-Tertikas \cite{FT}, %%@
Adimurthi-Grossi -Santra \cite{AGS},  Tertikas-Zographopoulos \cite{TZ},  and
Blanchet-Bonforte-Dolbeault-Grillo-Vasquez \cite{BBDGV}.

%Tertikas -Zographopoulos among others. The approach  also applies to improve the recent results of %%@Liskevich-Lyachova-Moroz on exterior domains and  will be developed in a forthcoming paper. 

 \end{abstract}

\section{Introduction} 

Ever since Br\'ezis-Vazquez \cite{BV} showed that Hardy's inequality can be improved  once %%@
restricted to a  smooth bounded domain $\Omega$ in $\R^n$, there was a flurry of activity about possible %%@
improvements  of the following type: 
\begin{equation}\label{gen-hardy.0}
\hbox{If $n\geq 3$\quad then \quad $\int_{\Omega}|\nabla u |^{2}dx - ( \frac{n-2}{2})^{2} %%@
\int_{\Omega}\frac{|u|^{2}}{|x|^{2}}dx\geq \int_{\Omega} V(x)|u|^{2}dx$ \quad for all $u \in %%@
H^{1}_{0}(\Omega)$,}
\end{equation} 
as well as its fourth order counterpart
\begin{equation}\label{gen-rellich.0}
\hbox{If $n\geq 5$\quad then \quad $\int_{\Omega}|\Delta u|^{2}dx -
\frac{n^2(n-4)^2}{16} \int_{\Omega}\frac{u^{2}}{|x|^{4}}dx\geq  
 \int_{\Omega} W(x)u^{2}dx$ \quad for $u  \in  H^{2}(\Omega)\cap H_{0}^{1}(\Omega)$},  
\end{equation}
where $V, W$ are certain explicit radially symmetric potentials of order lower than %%@
$\frac{1}{r^2}$ (for $V$) and  $\frac{1}{r^4}$ (for $W$).

In this paper, we provide an approach that completes, simplifies and improves most related %%@
results to-date regarding the Laplacian on Euclidean space as well as its powers. We also %%@
establish new inequalities some of which cover  critical dimensions such as $n=2$ for inequality %%@
(\ref{gen-hardy.0}) and $n=4$ for   (\ref{gen-rellich.0}).

\quad We start -- in section 2 -- by giving necessary and sufficient conditions  on positive %%@
radial functions  $V$ and $W$ on a ball $B$  in $R^{n}$,   so that the following inequality %%@
holds for some $c>0$:   
\begin{equation}\label{most.general.hardy}
\hbox{$\int_{B}V(x)|\nabla u |^{2}dx \geq c\int_{B} W(x)u^2dx$ for all $u \in %%@
C_{0}^{\infty}(B)$.}
\end{equation} 
Assuming that the ball $B$ has radius $R$ and that %%@
$\int^{R}_{0}\frac{1}{r^{n-1}V(r)}dr=+\infty$,  the condition is simply that the ordinary %%@
differential equation   
\begin{equation*}
\hbox{ $({\rm B}_{V,cW})$  \quad \quad \quad \quad \quad \quad \quad \quad \quad \quad \quad $y''(r)+(\frac{n-1}{r}+\frac{V_r(r)}{V(r)})y'(r)+\frac{cW(r)}{V(r)}y(r)=0$ \quad \quad \quad \quad \quad \quad \quad \quad \quad \quad \quad}
\end{equation*}
has a positive solution on the interval $(0, R)$. We shall call such a couple $(V, W)$  a 
{\it Bessel pair on $(0, R)$}. The {\it weight} of such a pair is then defined as
\begin{equation}
\hbox{$\beta (V, W; R)=\sup \big\{ c;\,  ({\rm B}_{V,cW})$ has a positive solution  on $(0, R)\big\} $.}
\end{equation}
This characterization makes an important connection between Hardy-type inequalities and the oscillatory behaviour of the above equations. For example, by 
using recent results on ordinary differential equations, we can then infer that an  integral %%@
condition on $V, W$ of the form
\begin{equation}
\limsup_{r\to 0}r^{2(n-1)}V(r)W(r)\big( \int^{R}_{r}\frac{d\tau}{\tau^{n-1}V(\tau)}\big)^2< %%@
\frac{1}{4}
\end{equation}
is sufficient (and  ``almost necessary") for $(V, W)$ to be a Bessel pair on a ball of sufficiently small radius $\rho$. 

\quad  Applied in particular, to a pair $(V, \frac{1}{r^2}V)$  where the function %%@
$\frac{rV'(r)}{V(r)}$ is assumed to decrease to $-\lambda$ on $(0, R)$, we obtain the following  %%@
extension of Hardy's inequality: 
 If $\lambda \leq n-2$, then  
\begin{equation}\label{v-hardy}
\hbox{$\int_{B}V(x)|\nabla u|^{2}dx\geq  %%@
(\frac{n-\lambda-2}{2})^2\int_{B}V(x)\frac{u^{2}}{|x|^2}dx$
\quad for all $u \in C^{\infty}_{0}(B)$}
\end{equation}
and $(\frac{n-\lambda-2}{2})^2$ is the best constant. The case where $V(x)\equiv 1$ is obviously the classical Hardy inequality and  when $V(x)=|x|^{-2a}$ for  $-\infty <a < \frac{n-2}{2}$, this  is a particular case  of the Caffarelli-Kohn-Nirenberg inequality. One can however apply the above criterium to obtain new inequalities such as the following: For $a,b> 0$
%, and $\alpha, \beta, m$ be real numbers. 
\begin{itemize} 
\item If $\alpha\beta>0$ and $m\leq \frac{n-2}{2}$, 
%or $\alpha, \beta<0$, and $2m-\alpha \beta %%@< n-2$, 
then  for all $u \in C^{\infty}_{0}(\R^n)$
\begin{equation}\label{GM-V1}
\int_{\R^n}\frac{(a+b|x|^{\alpha})^{\beta}}{|x|^{2m}}|\nabla u|^2dx\geq %%@
(\frac{n-2m-2}{2})^2\int_{\R^n}\frac{(a+b|x|^{\alpha})^{\beta}}{|x|^{2m+2}}u^2dx, 
\end{equation}
and $(\frac{n-2m-2}{2})^2$ is the best constant in the inequality.
\item If $\alpha \beta<0$ and $2m-\alpha \beta \leq n-2$, then for all $u \in C^{\infty}_{0}(\R^n)$
\begin{equation}\label{GM-V2}
\int_{\R^n}\frac{(a+b|x|^{\alpha})^{\beta}}{|x|^{2m}}|\nabla u|^2dx\geq (\frac{n-2m+\alpha %%@
\beta-2}{2})^2\int_{\R^n}\frac{(a+b|x|^{\alpha})^{\beta}}{|x|^{2m+2}}u^2dx, 
\end{equation}
and $(\frac{n-2m+\alpha %%@
\beta-2}{2})^2$ is the best constant in the inequality.
\end{itemize}
We can also extend some of the recent results of Blanchet-Bonforte-Dolbeault-Grillo-Vasquez \cite{BBDGV}.
% and answer some of their questions regarding best constants.  
%\begin{theorem}\label{GM-V} Let $a,b> 0$, and $\alpha, \beta$ be real numbers. 
\begin{itemize}
\item If $\alpha \beta <0$ and $-\alpha \beta \leq n-2$, %or $\alpha>0$, $\beta<0$, and $n\geq %%@2$,
 then for all $u \in C^{\infty}_{0}(\R^n)$
\begin{equation}\label{GM-V3}
\int_{\R^n}(a+b|x|^{\alpha})^{\beta}|\nabla u|^2dx\geq b^{\frac{2}{\alpha}}(\frac{n-\alpha %%@
\beta-2}{2})^2\int_{\R^n}(a+b|x|^{\alpha})^{\beta-\frac{2}{\alpha}}u^2dx,
\end{equation}
and $b^{\frac{2}{\alpha}}(\frac{n-\alpha 
\beta-2}{2})^2$ is the best constant in the inequality.
\item If $\alpha \beta >0$, 
%$\beta<0$, and $-\alpha \beta \leq n-2 $, or $\alpha>0$, $\beta>0$,
 and %%@
$n\geq 2$, then there exists a constant $C>0$ such that   for all $u %%@
\in C^{\infty}_{0}(\R^n)$
\begin{equation}\label{GM-V4}
\int_{\R^n}(a+b|x|^{\alpha})^{\beta}|\nabla u|^2dx\geq C %%@
\int_{\R^n}(a+b|x|^{\alpha})^{\beta-\frac{2}{\alpha}}u^2dx.
\end{equation}
Moreover, $b^{\frac{2}{\alpha}}(\frac{n-2}{2})^2\leq C\leq b^{\frac{2}{\alpha}}(\frac{n+\alpha %%@
\beta-2}{2})^2$.
\end{itemize}
On the other hand, by considering %%@
the  pair
\[
\hbox{$V(x)=|x|^{-2a}$\quad and \quad $W_{a,c}(x)= (\frac{n-2a-2}{2})^2|x|^{-2a-2} %%@
+c|x|^{-2a}W(x)$}
\]
we get the following improvement of the Caffarelli-Kohn-Nirenberg inequalities:
\begin{equation}\label{CKN}
 \int_{B}|x|^{-2a}|\nabla u |^{2}dx -  (\frac{n-2a-2}{2})^2\int_{B}|x|^{-2a-2}u^2 dx\geq %%@
c\int_{B}|x|^{-2a} W(x)u^2dx \quad \hbox{for all $u \in C_{0}^{\infty}(B)$}
\end{equation} 
if and only if  the following ODE 
\begin{equation*}\label{ODE111}
\hbox{ $({\rm B}_{cW})$  \quad \quad \quad \quad\quad \quad \quad \quad \quad \quad \quad \quad \quad \quad \quad $y''+\frac{1}{r}y'+c W(r)y=0$ \quad \quad \quad \quad \quad \quad \quad \quad \quad \quad \quad\quad \quad \quad \quad}
\end{equation*}
 has a positive solution on $(0, R)$. Such a function $W$ will be called a {\it Bessel potential} on $(0, R)$. This type of  characterization was established recently by the authors \cite{GM1} in the case  
where $a=0$, yielding in particular  the recent improvements of Hardy's inequalities (on bounded domains) established by %%@
Brezis-V\'{a}zquez \cite{BV},  Adimurthi et al.   \cite{ACR}, and Filippas-Tertikas \cite{FT}.  %%@
Our results here include in addition those proved by Wang-Willem \cite{WW} in the case where $a< %%@
\frac{n-2}{2}$ and $W(r)=\frac{1}{r^2(\ln\frac{R}{r})^2}$, but also cover the  previously unknown limiting case corresponding to $a=  
\frac{n-2}{2}$ as well as the critical dimension  $n=2$.

More importantly, we establish here  that Bessel pairs lead to a myriad of optimal Hardy-Rellich %%@
inequalities of arbitrary high order, therefore extending and completing a series of new results %%@
by Adimurthi et al. \cite{AGS}, Tertikas-Zographopoulos  \cite{TZ} and others. They are mostly %%@
based on the following theorem which summarizes the main thrust of this paper.

\begin{theorem}  Let $V$ and $W$ be positive radial $C^1$-functions   on $B\backslash \{0\}$, %%@
where $B$ is a ball centered at zero with radius $R$ in $\R^n$ ($n \geq 1$) such that  %%@
$\int^{R}_{0}\frac{1}{r^{n-1}V(r)}dr=+\infty$ and $\int^{R}_{0}r^{n-1}V(r)dr<+\infty$. The %%@
following statements are then equivalent:

\begin{enumerate}

\item $(V, W)$ is a Bessel pair on $(0, R)$ and $\beta (V, W; R) \geq 1$. 

\item $ \int_{B}V(x)|\nabla u |^{2}dx \geq \int_{B} W(x)u^2dx$ for all $u \in C_{0}^{\infty}(B)$.

\item If $\lim_{r \rightarrow 0}r^{\alpha}V(r)=0$ for some $\alpha< n-2$, then the above are %%@
equivalent to
\[  
\hbox{$\int_{B}V(x)|\Delta u |^{2}dx \geq  \int_{B} W(x)|\nabla  %%@
u|^{2}dx+(n-1)\int_{B}(\frac{V(x)}{|x|^2}-\frac{V_r(|x|)}{|x|})|\nabla u|^2dx$ \quad  for all %%@
radial $u \in  
C^{\infty}_{0,r}(B)$.}
\]
 
 \item If in addition, $W(r)-\frac{2V(r)}{r^2}+\frac{2V_r(r)}{r}-V_{rr}(r)\geq 0$ on $(0, R)$, %%@
then the above are equivalent to
\[  
\hbox{$\int_{B}V(x)|\Delta u |^{2}dx \geq  \int_{B} W(x)|\nabla  %%@
u|^{2}dx+(n-1)\int_{B}(\frac{V(x)}{|x|^2}-\frac{V_r(|x|)}{|x|})|\nabla u|^2dx$ \quad  for all  $u %%@
\in  
C^{\infty}_{0}(B)$.}
\]
\end{enumerate}
\end{theorem} 
In other words, one can obtain as many Hardy and Hardy-Rellich type inequalities as one can %%@
construct  Bessel pairs on $(0, R)$.  The relevance of the above result stems from the fact that %%@
there are plenty of such pairs that are easily identifiable. Indeed, even the class of {\it Bessel %%@
potentials} --equivalently those $W$ such that $\left(1, (\frac{n-2}{2})^2|x|^{-2} +cW(x)\right)$ is a Bessel pair-- is quite  rich and contains  %%@
several important potentials.  Here are  some of the  most relevant properties --to be established %%@
in an appendix-- of the class  of $C^1$ Bessel potentials  $W$ on $(0, R)$, that we shall denote %%@
by ${\cal B}(0, R)$.

First, the class is a closed convex {\it solid} subset of $C^1(0, R)$, that is 
if $W\in {\cal B}(0, R)$ and $0\leq V\leq W$, then $V\in {\cal B}(0, R)$. The "weight" of each %%@
$W\in {\cal B}(R)$, that is 
\begin{equation}
\hbox{$\beta (W; R)=\sup\big\{c>0;\, (B_{cW})$ has a positive solution on $(0, R)\big\},$}
\end{equation}
 will be an important ingredient for computing the best constants in corresponding functional %%@
inequalities. Here are some basic examples of Bessel potentials and their corresponding weights.

\begin{itemize}
\item $  W \equiv 0$ is a Bessel potential on $(0,  R)$ for any $R>0$. 
\item $  W \equiv 1$ is a Bessel potential on $(0,  R)$ for any $R>0$, and $\beta (1; %%@
R)=\frac{z_0^2}{R^2}$  where $z_{0}=2.4048...$ is the first zero of the Bessel function $J_0$.

\item   If $a<2$, then there exists $R_a>0$ such that $W (r)=r^{-a}$ is  a Bessel potential on $( %%@
0, R_a)$. 

\item  For  $k\geq 1$, $R>0$ and $\rho=R( e^{e^{e^{.^{.^{e((k-1)-times)}}}}} )$, let  $ W_{k, %%@
\rho} (r)=\Sigma_{j=1}^k\frac{1}{r^{2}}\big(\prod^{j}_{i=1}log^{(i)}\frac{\rho}{r}\big)^{-2}$ %%@
where the functions $log^{(i)}$ are defined  iteratively  as follows:  $log^{(1)}(.)=log(.)$ and  %%@
for $k\geq 2$,  $log^{(k)}(.)=log(log^{(k-1)}(.))$.  $ W_{k, \rho}$ is then a Bessel potential on %%@
$(0, R)$ with $\beta (W_{k, \rho}; R)=\frac{1}{4}$. 
 
\item  For $k\geq 1$, $R>0$ and $\rho\geq R$, define $\tilde W_{k; \rho} %%@
(r)=\Sigma_{j=1}^k\frac{1}{r^{2}}X^{2}_{1}(\frac{r}{\rho})X^{2}_{2}(\frac{r}{\rho}) \ldots  
X^{2}_{j-1}(\frac{r}{\rho})X^{2}_{j}(\frac{r}{\rho})$ where   the functions $X_i$ are defined %%@
iteratively  as follows:
 $X_{1}(t)=(1-\log(t))^{-1}$ and for $k\geq 2$, $ X_{k}(t)=X_{1}(X_{k-1}(t))$. Then again $ \tilde %%@
W_{k, \rho}$ is a Bessel potential on $(0, R)$ with $\beta (\tilde W_{k, \rho}; R)=\frac{1}{4}$. 
 
 \item More generally, if $W$ is any positive  function on $\R$ such that 
$\liminf\limits_{r\rightarrow 0} \ln(r)\int^{r}_{0} sW(s)ds>-\infty$,
then for every $R>0$, there exists $\alpha:=\alpha(R)>0$ such that $W_\alpha(x):=\alpha^2W(\alpha %%@
x)$ is a Bessel potential on $(0, R)$. 
\end{itemize}

What is remarkable is that the class of Bessel potentials $W$ is also the one that leads to %%@
optimal improvements for fourth order inequalities (in dimension $n\geq 3$) of the following type: 

\begin{equation}\label{gen-second.hardy}
\hbox{$ \int_{B}|\Delta u  
|^{2}dx - C(n) \int_{B}\frac{|\nabla u|^{2}}{|x|^{2}}dx\geq  c(W, R)\int_{B}  
W(x)|\nabla u|^{2}dx$ \quad for all $u \in H^{2}_{0}(B)$,}  
  \end{equation}
  where $C(3)=\frac{25}{36}$, $C(4)=3$ and $C(n)=\frac{n^2}{4}$ for $n\geq 5$. 
The case when $W\equiv \tilde W_{k, \rho}$ and $n\geq 5$ was recently  established  by %%@
Tertikas-Zographopoulos \cite{TZ}. Note that  $W$ can be chosen to be any one of the examples of %%@
Bessel potentials listed above. Moreover, both $C(n)$ and the weight $\beta (W; R)$ are the best %%@
constants in the above inequality. 

 Appropriate combinations of  (\ref{most.general.hardy}) and (\ref{gen-second.hardy}) then lead  %%@
to  a myriad of Hardy-Rellich inequalities in dimension $n\geq 4$. For example, if $W$ is a  %%@
Bessel potential on $(0, R)$ such that the function $r\frac{W_r(r)}{W(r)}$ decreases to %%@
$-\lambda$, and if $\lambda \leq n-2$, then we have for all $u \in C^{\infty}_{0}(B_{R})$
  \begin{equation} \label{Rellich.1}
\int_{B}|\Delta u|^{2}dx -
\frac{n^2(n-4)^2}{16}\int_{B}\frac{u^2}{|x|^4}dx\geq %%@
\big(\frac{n^2}{4}+\frac{(n-\lambda-2)^2}{4}\big)
 \beta (W; R)\int_{B}\frac{W(x)}{|x|^2}u^2 dx.
\end{equation}

By applying (\ref{Rellich.1}) to the various examples of Bessel functions listed above,  one  %%@
improves in many ways the recent results of Adimurthi et al. \cite{AGS} and those by %%@
Tertikas-Zographopoulos  \cite{TZ}. Moreover, besides covering the critical dimension $n=4$, we %%@
also establish that the best constant is $(1+\frac{n(n-4)}{8})$ for all the potentials $W_k$ and %%@
$\tilde W_k$ defined above. For example we have for $n\geq 4$, 
\begin{equation}
 \int_{B}|\Delta u(x) |^{2}dx \geq \frac{n^2(n-4)^2}{16}\int_{B}\frac{u^2}{|x|^4} %%@
dx+(1+\frac{n(n-4)}{8})\sum^{k}_{j=1}\int_{B}\frac{u^2}{|x|^4}\big( %%@
\prod^{j}_{i=1}log^{(i)}\frac{\rho}{|x|}\big)^{-2}dx.
 \end{equation}
% and
% \begin{equation}
% \int_{B}|\Delta u(x) |^{2}dx \geq \frac{n^2(n-4)^2}{16}\int_{B}\frac{u^2}{|x|^4} %%@
%dx+(1+\frac{n(n-4)}{8})\sum^{k}_{j=1}\int_{B}\frac{u^2}{|x|^4}X^{2}_{1}(\frac{|x|}{\rho})X^{2%%@
%%@
%}_{2}(\frac{r}{\rho}) \ldots  
%X^{2}_{j-1}(\frac{|x|}{\rho})X^{2}_{j}(\frac{|x|}{\rho})dx.
% \end{equation}
 
More generally, we show that for any  $m<\frac{n-2}{2}$, and any $W$ Bessel potential on a ball %%@
$B_{R}\subset R^n$ of radius $R$, the following inequality holds for all $u \in C^{\infty}_{0}(B_{R})$
\begin{equation}\label{gm-hr.00}
\int_{B_{R}}\frac{|\Delta u|^2}{|x|^{2m}}\geq a_{n,m}\int_{B_{R}}\frac{|\nabla %%@
u|^2}{|x|^{2m+2}}dx+\beta(W; R)\int_{B_{R}}W(x)\frac{|\nabla u|^2}{|x|^{2m}}dx,
\end{equation}
where $a_{m,n}$ and $\beta(W; R)$  are best constants that we compute in the appendices for all %%@
$m$ and $n$ and  for many Bessel potentials $W$. Worth noting is Corollary \ref{radial} where we %%@
show that inequality (\ref{gm-hr.00}) restricted to radial functions in  $C^{\infty}_{0}(B_{R})$ %%@
holds with a best constant equal to $(\frac{n+2m}{2})^2$, but that $a_{n,m}$ can however be %%@
strictly smaller than $(\frac{n+2m}{2})^2$ in the non-radial case. These results improve %%@
considerably Theorem 1.7, Theorem 1.8, and Theorem 6.4 in \cite{TZ}.

 We also establish a more general version of equation (\ref{Rellich.1}). Assuming again that  
$\frac{rW'(r)}{W(r)}$  decreases to $-\lambda$ on $(0, R)$, and provided  $m\leq \frac{n-4}{2}$ %%@
and $\lambda \leq n-2m-2$, we then have  for all $u \in C^{\infty}_{0}(B_R)$, 
\begin{eqnarray}\label{ex-gen-hr}
\int_{B_R}\frac{|\Delta u|^{2}}{|x|^{2m}}dx &\geq& \beta_{n,m}\int_{B_R}\frac{u^2}{|x|^{2m+4}}dx %%@
\nonumber\\
&&\quad+\beta (W;  R)(\frac{(n+2m)^2}{4}+\frac{(n-2
m-\lambda-2)^2}{4})
\int_{B_R}\frac{W(x)}{|x|^{2m+2}}u^2 dx,
\end{eqnarray}
where again the best constants $\beta_{n,m}$ are computed in section 3.  This completes the %%@
results in Theorem 1.6 of \cite{TZ}, where the inequality is established for $n\geq 5$, $0\leq m < \frac{n-4}{2}$, and the particular potential ${\tilde W}_{k,\rho}$.  

Another inequality that relates the Hessian integral to the Dirichlet energy is the %%@
following: 
Assuming $-1<m\leq \frac{n-4}{2}$ and $W$ is a Bessel potential on a ball $B$ of radius $R$ in %%@
$R^n$, then  for all $u \in C^{\infty}_{0}(B)$,
\begin{eqnarray}\label{hrs}
\int_{B}\frac{|\Delta u|^{2}}{|x|^{2m}}dx- %%@
\frac{(n+2m)^{2}(n-2m-4)^2}{16}\int_{B}\frac{u^2}{|x|^{2m+4}}dx&\geq&
\beta (W; R)\frac{(n+2m)^2}{4}
\int_{B}\frac{W(x)}{|x|^{2m+2}}u^2 dx \quad \quad \quad \nonumber \\
&&+ \beta(|x|^{2m}; R)||u||_{H^1_0}.
\end{eqnarray}
This improves considerably Theorem A.2. in \cite{AGS} where it is established -- for  $m=0$ and without best constants -- with the potential  $W_{1, \rho}$
%=\frac{1}{r^2(ln\frac{R}{r})^2}$ 
in dimension $n\geq 5$,  and the potential $W_{2, \rho}$ when $n=4$.
%=\frac{1}{r^2(ln\frac{R}{r})^2}$

Finally, we establish several higher order Rellich inequalities for integrals of the form $\int_{B_{R}}\frac{|\Delta^{m}u|^2}{|x|^{2k}}dx$, improving in many ways several recent results in \cite{TZ}. 

% a sample of which is the %%@
%following: If  $m \in N$, $1\leq l\leq m$, and $2k+4m\leq n$, then  for all $u \in %%@
%C^{\infty}_{0}(B)$
%\begin{eqnarray}\label{nuts.0}
%\int_{B_{R}}\frac{|\Delta^{m}u|^2}{|x|^{2k}}dx &\geq& %%@
%\prod\limits_{i=1}^{l}\frac{a_{_{n,k+2i-2}}(n-2k-4i)^2}{4}\int_{B_{R}}\frac{|\Delta^{m-l}u|^2}{|x|%%@
%^{2k+4l}}dx\\
%&+&\beta(W; %%@
%R)\sum\limits_{i=1}^{l}\prod\limits_{j=1}^{l-1}\frac{a_{_{n,k+2j-2}}(n-2k-4j)^2}{4}\int_{B_{R
%}}W(x)\frac{|\nabla \Delta^{m-i}u|^2}{|x|^{2k+4i-4}}dx\nonumber\\
%&+&\beta(W; %%@
%R)\sum\limits_{i=1}^{l}a_{_{n,k+2i-2}}\prod\limits_{j=1}^{l-1}\frac{a_{_{n,k+2j-2}}(n-2k-4j)^2}
%{4})\int_{B_{R}}W(x)\frac{|\Delta^{m-i}u|^2}{|x|^{2k+4i-2}}dx,\nonumber
%\end{eqnarray}
%where $a_{n,m}$ are the best constants in inequality (\ref{gm-hr.00}). This improves Theorem 1.10 %%@
%in \cite{TZ} where it is proved for $k=0$, $l\leq  \frac{-n+8+2\sqrt{n^2-n+1}}{12}$, $4m<n$, and $W={\tilde W}_{k, \rho}$. %%@
%Note that even for $k=0$, inequality (\ref{nuts.0}) shows that  we   can drop the first condition %%@
%and replace the second one by $4m\leq n$  in  Theorem 1.10 of \cite{TZ}. 

The approach  can also be used to improve the recent results of  
Liskevich-Lyachova-Moroz \cite{LLM} on exterior domains and  will be developed in a forthcoming paper. 

\section{General Hardy Inequalities}

Here is the main result of this section.
\begin{theorem} \label{main} Let $V$ and $W$ be positive radial $C^1$-functions on  %%@
$B_R\backslash \{0\}$, where $B_R$ is a ball centered at %%@
zero with radius $R$ ($0<T\leq +\infty$) in $\R^n$ ($n \geq 1$). Assume that $\int^{a}_{0}\frac{1}{r^{n-1}V(r)}dr=+\infty$ and %%@
$\int_{0}^{a}r^{n-1}V(r)dr<\infty$ for some $0<a<R$. Then the following two statements are equivalent:

\begin{enumerate}
\item The ordinary differential equation 
\[
\hbox{ $({\rm B}_{V,W})$  \quad \quad \quad \quad \quad  \quad \quad \quad \quad \quad \quad \quad %%@
\quad \quad \quad $y''(r)+(\frac{n-1}{r}+\frac{V_r(r)}{V(r)})y'(r)+\frac{W(r)}{V(r)}y(r)=0$  \quad %%@
\quad \quad \quad \quad  \quad \quad \quad \quad \quad  \quad \quad \quad \quad \quad \quad \quad %%@
\quad \quad \quad}
\]
 has a positive solution on the interval $(0, R]$ (possibly with $\varphi(R)=0)$.
 
\item For all $u \in C_{0}^{\infty}(B_R)$
\begin{equation*}\label{2dim-in}
\hbox{ $({\rm H}_{V,W})$ \quad \quad \quad \quad \quad \quad \quad \quad \quad \quad \quad \quad %%@
\quad   $\int_{B_R}V(x)|\nabla u(x) |^{2}dx \geq \int_{B_R} W(x)u^2dx$.\quad \quad \quad \quad \quad %%@
\quad \quad  \quad \quad \quad \quad \quad \quad \quad \quad \quad \quad}
\end{equation*} 
\end{enumerate}

\end{theorem} 

%\begin{theorem} \label{main-Rn} Let $V$ and $W$ be positive radial $C^1$-functions on a ball %%@
%$\R^n\backslash \{0\}$, with $n \geq 1$. Assume that $\int^{a}_{0}\frac{1}{r^{n-1}V(r)}dr=+\infty$ %%@
%and $\int_{0}^{a}r^{n-1}V(r)dr<\infty$ for some $a>0$, then the following two statements are %%@
%equivalent:
%\begin{enumerate}
%\item The ordinary differential equation 
%\[
%\hbox{ $({\rm B'}_{V,W})$  \quad \quad \quad \quad \quad  \quad \quad \quad \quad \quad \quad %%@
%\quad \quad \quad \quad $y''(r)+(\frac{n-1}{r}+\frac{V_r(r)}{V(r)})y'(r)+\frac{W(r)}{V(r)}y(r)=0$  %%@
%\quad \quad \quad \quad \quad  \quad \quad \quad \quad \quad  \quad \quad \quad \quad \quad \quad %%@
%\quad \quad \quad \quad}
%\]
% has a positive solution on the interval $(0, \infty)$.
 %\item For all $u \in C_{0}^{\infty}(\R^n)$
%\begin{equation*}\label{2dim-in}
%\hbox{ $({\rm H'}_{V,W})$ \quad \quad \quad \quad \quad \quad \quad \quad \quad \quad \quad \quad %%@
%\quad   $\int_{\R^n}V(x)|\nabla u(x) |^{2}dx \geq \int_{\R^n} W(x)u^2dx$.\quad \quad \quad \quad %%@
%\quad \quad \quad  \quad \quad \quad \quad \quad \quad \quad \quad \quad \quad}
%\end{equation*} 
%\end{enumerate}
%\end{theorem} 
Before proceeding with the proofs, we note the following immediate but useful corollary.

\begin{corollary} Let $V$ and $W$ be positive radial $C^1$-functions on $B\backslash \{0\}$,  %%@
where $B$ is a ball with radius $R$ in $\R^n$ ($n \geq 1$) and centered at zero,  such that %%@
$\int^{R}_{0}\frac{1}{r^{n-1}V(r)}dr=+\infty$ and $\int^{R}_{0}r^{n-1}V(r)dr<\infty$. Then $(V, %%@
W)$ is a Bessel pair on $(0, R)$ if and only if for all $u \in C_{0}^{\infty}(B_R)$, we have
 \[
 \int_{B_R}V(x)|\nabla u |^{2}dx \geq \beta (V, W; R) \int_{B_R} W(x)u^2dx, 
 \]
with $\beta (V, W; R)$ being the best constant.
\end{corollary}

For the proof of Theorem \ref{main}, we shall need the following lemmas.

\begin{lemma}\label{main-lem} Let $\Omega$ be a smooth bounded domain in $R^n$ with $n\geq 1$ and %%@
let $\varphi \in C^{1}(0,R:=\sup_{x \in \partial \Omega}|x|)$ be a positive solution of the %%@
ordinary differential equation 
\begin{equation}\label{ODEv}
y''+(\frac{n-1}{r}+\frac{V_r(r)}{V(r)})y'+\frac{W(r)}{V(r)}y=0,
\end{equation}
on $(0,R)$ for some $V(r),W(r)\geq 0$ where $\int^{R}_{0}\frac{1}{r^{n-1}V(r)}dr=+\infty$ and %%@
$\int^{R}_{0}r^{n-1}V(r)dr<\infty$. Setting $\psi(x)=\frac{u(x)}{\varphi(|x|)}$ for any $u \in %%@
C^{\infty}_{0}(\Omega)$, we then have the following properties:
\begin{enumerate}
\item $\int_{0}^{R}r^{n-1}V(r)(\frac{\varphi'(r)}{\varphi(r)})^2dr<\infty$ and $\lim_{r\rightarrow %%@
0}r^{n-1}V(r)\frac{\varphi'(r)}{\varphi(r)}=0.$
\item $\int_{\Omega}V(|x|)(\varphi'(|x|))^2\psi^2(x)dx<\infty.$
\item $\int_{\Omega}V(|x|)\varphi^2(|x|)|\nabla \psi|^2(x)dx<\infty.$
\item $|\int_{\Omega}V(|x|)\varphi'(|x|)\varphi(|x|)\psi(x)\frac{x}{|x|}.\nabla %%@
\psi(x)dx|<\infty$.
\item $\lim_{r\rightarrow 0}|\int_{\partial  
B_{r}}V(|x|)\varphi'(|x|)\varphi(|x|)\psi^2(x)ds|=0,$ where 
$B_{r}\subset \Omega$ is a ball of radius $r$ centered at $0$. 
\end{enumerate}
\end{lemma}
\noindent {\bf Proof:}  $1)$ Setting 
$x(r)=r^{n-1}V(r)\frac{\varphi'(r)}{\varphi(r)}$, we have
\[ %%@
r^{n-1}V(r)x'(r)+x^2(r)=\frac{r^{2(n-1)}V^{2}(r)}{\varphi}(\varphi''(r)+(\frac{n-1}{r}+\frac{V_r(r
)}{V(r)})\varphi'(r))=-
\frac{r^{2(n-1)}V(r)W(r)}{\varphi(r)}\leq 0, \ \ \ \ 0<r<R.\]
Dividing by $r^{n-1}V(r)$ and integrating once, we obtain
\begin{equation}\label{lem-eq}
x(r)\geq \int_{r}^{R}\frac{|x(s)|^2}{s^{n-1}V(s)}ds+x(R).
\end{equation}
To prove that $\lim_{r\rightarrow 0}G(r)<\infty$, where $G(r):=\int_{r}^{R}\frac{x^{2}(s)}{s^{n-1}V(s)}ds$,  we assume the contrary and 
%  Assuming that  $G(r):=\int_{r}^{R}\frac{x^{2}(s)}{s^{n-1}V(s)}ds \rightarrow \infty$ 
%as $r\rightarrow 0$, we get from 
use (\ref{lem-eq}) to write that
\begin{equation*}
(-r^{n-1}V(r))G'(r))^{\frac{1}{2}}\geq G(r)+x(R).
\end{equation*}
%Note that $G$ goes to infinity as $r$ goes to zero. 
Thus, for $r$ sufficiently small we have
$-r^{n-1}V(r)G'(r)\geq \frac{1}{2}G^{2}(r)$
and hence,
$(\frac{1}{G(r)})'\geq\frac{1}{2r^{n-1}V(r)}$, which contradicts the fact that $G(r)$ goes to %%@
infinity as $r$ tends to zero. 

Also in view of (\ref{lem-eq}), we have that  $x_{0}:=\lim_{r\rightarrow %%@
0}x(r)$ exists, and since $\lim_{r\rightarrow 0}G(r)<\infty$, we necessarily have $x_{0}=0$ and 1) is proved. 

For assertion 2), we  use  $1)$ to see that
\[\int_{\Omega}V(|x|)(\varphi'(|x|))^2\psi^2(x)dx \leq ||u||^{2}_{\infty}\int_{\Omega}V(|x|) %%@
\frac{(\varphi'(|x|))^2}{\varphi^2(|x|)}dx<\infty.\]

3)\,  Note that 
\[\hbox{$|\nabla \psi(x)|\leq \frac{|\nabla %%@
u(x)|}{\varphi(|x|)}+|u(x)|\frac{|\varphi'(|x|)|}{\varphi^2(|x|)}\leq %%@
\frac{C_{1}}{\varphi(|x|)}+C_{2}\frac{|\varphi'(|x|)|}{\varphi^2(|x|)}$, \ \ for all $x \in \Omega %%@
$},\]
where $C_{1}=\max_{x \in \Omega}|\nabla u|$ and $C_{2}=\max_{x \in \Omega}|u|$.
Hence we have
\begin{eqnarray*}
\int_{\Omega}V(|x|)\varphi^2(|x|)|\nabla \psi|^2(x)dx&\leq& \int_{\omega} %%@
V(|x|)\frac{(C_{1}\varphi(|x|)+C_{2}\varphi'(|x|))^2}{\varphi^2(|x|)}dx\\
&=&\int_{\Omega}C^{2}_{1}V(|x|)dx+\int_{\Omega}2C_{1}C_{2}\frac{|\varphi'(|x|)|}{\varphi(|x|)}
V(|x|)dx+\int_{\Omega}C^{2}_{2}(\frac{\varphi'(|x|)}{\varphi(|x|)})^2V(|x|) dx\\
&\leq& %%@
L_{1}+2C_{1}C_{2}\big(\int_{\Omega}V(|x|)(\frac{\varphi'(|x|)}{\varphi(|x|)})^2dx\big)^{\frac{1
}{2}} \big(\int_{\Omega}V(|x|)dx\big)^{\frac{1}{2}}+L_{2}\\
&<& \infty,
\end{eqnarray*}
which proves $3)$. 

$4)$ now follows from $2)$ and $3)$ since
\[ V(|x|)|\nabla u|^2= V(|x|)(\varphi'(|x|))^2 %%@
\psi^2(x)+2V(|x|)\varphi'(|x|)\varphi(|x|)\psi(x)\frac{x}{|x|}.\nabla %%@
\psi(x)+V(|x|)\varphi^2(|x|)|\nabla \psi|^2.
\]
Finally,  $5)$ follows from $1)$ since
\begin{eqnarray*}
|\int_{\partial  
B_{r}}V(|x|)\varphi'(|x|)\varphi(|x|)\psi^2(x)ds|&<&||u||^{2}_{\infty}|\int_{\partial %%@
B_{r}}V(|x|)\frac{\varphi'(|x|)}{\varphi(|x|)}ds\\
&=&||u||^{2}_{\infty}V(r)\frac{|\varphi'(r)|}{\varphi(r)}\int_{\partial B_{r}}1 ds\\
&=&n\omega_{n}||u||^{2}_{\infty}r^{n-1}V(r)\frac{|\varphi'(r)|}{\varphi(r)}.
\end{eqnarray*}
%The proof of the lemma is complete. \hfill $\Box$

\begin{lemma} \label{super} Let $V$ and $W$ be positive radial $C^1$-functions  on a ball %%@
$B\backslash \{0\}$, where $B$ is a ball with radius $R$ in $\R^n$ ($n \geq 1$) and centered at %%@
zero. Assuming
\begin{eqnarray*}
\hbox{$\int_{B}\left(V(x)|\nabla u|^{2}-W(x)|u|^{2}\right)dx\geq 0$ for all $u \in %%@
C_{0}^{\infty}(B)$, }
\end{eqnarray*}
 then there exists a $C^{2}$-supersolution to the following linear elliptic equation
\begin{eqnarray}\label{pde}
-{\rm div}(V(x)\nabla u)-W(x)u&=&0, \ \ \ \ {\rm in} \ \  B, \\
u&>&0 \ \ \quad {\rm in} \ \  B \setminus \{0\}, \\
 u&=&0 \ \quad {\rm in} \ \  \partial B. 
\end{eqnarray}
 \end{lemma}
{\bf Proof:} Define
 \begin{eqnarray*}
\lambda_{1}(V):=\inf \{\frac{\int_{B}V(x)|\nabla\psi|^{2}- %%@
W(x)|\psi|^{2}}{\int_{B}|\psi|^{2}}; \ \ \psi \in C^{\infty}_{0}(B \setminus \{0\}) \}.
\end{eqnarray*}
By our assumption $\lambda_1(V)\geq 0$. Let $(\phi_{n}, \lambda^{n}_{1})$ be the first eigenpair %%@
for the problem
\begin{eqnarray*}
(L-\lambda_{1}(V)-\lambda^{n}_{1})\phi_{n}&=&0 \ \ on \ \ B \setminus B_{\frac{R}{n}}\\
\phi_n&=&0 \ \ on\ \ \partial (B \setminus B_{\frac{R}{n}}),
\end{eqnarray*}
where $Lu=-{\rm div}(V(x)\nabla u)- W(x) u$, and $B_{\frac{R}{n}}$ is a ball of radius %%@
$\frac{R}{n}$, $n\geq 2$ . The eigenfunctions can be chosen in such a way that $\phi_{n}>0$ on %%@
$B \setminus B_{\frac{R}{n}}$ and $\varphi_{n}(b)=1$, for some $b \in B$ with %%@
$\frac{R}{2}<|b|<R$.

Note that $\lambda^{n}_{1}\downarrow 0$ as $n \rightarrow \infty$. Harnak's inequality yields that %%@
for any compact subset $K$, $\frac{{\rm max}_{K}\phi_{n}}{{\rm min}_{K}\phi_{n}}\leq C(K)$ with %%@
the later constant being independant of $\phi_{n}$.  Also standard elliptic estimates also yields %%@
that the family $(\phi_{n})$ have also uniformly bounded derivatives on the compact sets %%@
$B-B_{\frac{R}{n}}$.  \\
Therefore, there exists a subsequence $(\varphi_{n_{l_{2}}})_{l_{2}}$ of ($\varphi_{n})_{n}$ such %%@
that $(\varphi_{n_{l_{2}}})_{l_{2}}$ converges to some $\varphi_{2} \in C^{2}(B \setminus %%@
B(\frac{R}{2}))$.  Now consider  $(\varphi_{n_{l_{2}}})_{l_{2}}$ on $B \setminus %%@
B(\frac{R}{3})$.  Again there exists a subsequence  $(\varphi_{n_{l_{3}}})_{l_{3}}$ of %%@
$(\varphi_{n_{l_{2}}})_{l_{2}}$ which converges to $\varphi_{3} \in C^{2}( B \setminus %%@
B(\frac{R}{3}))$,  and $\varphi_{3}(x)=\varphi_{2}(x)$ for all $x \in B \setminus %%@
B(\frac{R}{2})$. By repeating this argument we get a supersolution $\varphi \in C^{2}( B %%@
\setminus\{ 0\})$ i.e. $L\varphi \geq 0$, such that $\varphi>0$ on $B \setminus \{0\}$. %%@
\hfill $\square$\\
 
\noindent {\bf Proof of Theorem \ref{main}:} First we prove that 1) implies 2). Let $\phi \in %%@
C^{1}(0,R]$ be a solution of $(B_{V,W})$   such that $\phi (x)>0$ for all $x \in (0,R)$. Define %%@
$\frac{u(x)}{\varphi(|x|)}= \psi(x)$. Then
\[ |\nabla u|^2=(\varphi'(|x|))^2 \psi^2(x)+2\varphi'(|x|)\varphi(|x|)\psi(x)\frac{x}{|x|}.\nabla %%@
\psi+\varphi^2(|x|)|\nabla \psi|^2.\]
Hence,
\[ V(|x|)|\nabla u|^2\geq V(|x|)(\varphi'(|x|))^2 %%@
\psi^2(x)+2V(|x|)\varphi'(|x|)\varphi(|x|)\psi(x)\frac{x}{|x|}.\nabla \psi(x).\]
Thus, we have
\[\int_{B} V(|x|)|\nabla u|^2 dx \geq \int_{B} V(|x|) (\varphi'(|x|))^2 \psi^2(x)dx+ \int_{B}2 %%@
V(|x|) \varphi'(|x|)\varphi(|x|)\psi(x) \frac{x}{|x|}.\nabla \psi dx.\]
Let $B_{\epsilon}$ be a ball of radius $\epsilon$ centered at the origin. Integrate by parts to %%@
get
\begin{eqnarray*}
\int_{B} V(|x|)|\nabla u|^2 dx &\geq& \int_{B} V(|x|) (\varphi'(|x|))^2 %%@
\psi^2(x)dx+\int_{B_{\epsilon}}2V(|x|) \varphi'(|x|)\varphi(|x|)\psi(x) \frac{x}{|x|}.\nabla \psi %%@
dx\\
&+&\int_{B\backslash B_{\epsilon}}2 V(|x|) \varphi'(|x|)\varphi(|x|)\psi(x) \frac{x}{|x|}.\nabla %%@
\psi dx\\
&=&\int_{B_{\epsilon}} V(|x|) (\varphi'(|x|))^2 \psi^2(x)dx+\int_{B_{\epsilon}}2V(|x|) %%@
\varphi'(|x|)\varphi(|x|)\psi(x) \frac{x}{|x|}.\nabla \psi dx\\
&-&\int_{B\backslash B_{\epsilon}}\left\{\big(V(|x|) %%@
\varphi''(|x|)\varphi(|x|)+(\frac{(n-1)V(|x|)}{r}+V_r(|x|))\varphi'(|x|)\varphi(|x|)\big)\psi^2(x)
\right\}dx\\
&+&\int_{\partial (B \backslash B_{\epsilon})}V(|x|)\varphi'(|x|)\varphi(|x|)\psi^2(x)ds
\end{eqnarray*}
Let $\epsilon \rightarrow 0$ and use Lemma $\ref{main-lem}$ and the fact that $\phi$ is a solution %%@
of   $(D_{v,w})$ to get
\begin{eqnarray*}
\int_{B} V(|x|)|\nabla u|^2 dx &\geq& %%@
-\int_{B}[V(|x|)\varphi''(|x|)+(\frac{(n-1)V(|x|)}{r}+V_r(|x|))\varphi'(|x|)]\frac{u^2(x)}{
\varphi(|x|)}dx\\
&=&\int_{B}W(|x|)u^2(x)dx.
\end{eqnarray*}
To show that 2) implies 1), we assume that inequality (${\rm H}_{V,W}$) holds on a ball $B$ of %%@
radius $R$, and then apply Lemma \ref{super} to obtain a $C^{2}$-supersolution for the equation %%@
(\ref{pde}). Now take the surface average of $u$, that is
\begin{equation}\label{sup-eq}
y(r)=\frac{1}{n\omega_{w} r^{n-1}}\int_{\partial B_{r}} u(x)dS=\frac{1}{n\omega_{n}} %%@
\int_{|\omega|=1}u(r\omega)d\omega >0,
\end{equation}
where $\omega_{n}$ denotes the volume of the unit ball in $R^{n}$. We may assume that the unit %%@
ball is contained in $B$ (otherwise we just use a smaller ball).  We clearly have
\begin{equation}
y''(r)+\frac{n-1}{r}y'(r)= \frac{1}{n\omega_{n}r^{n-1}}\int_{\partial B_{r}}\Delta u(x)dS.
\end{equation}
Since $u(x)$  is a supersolution of (\ref{pde}), we have
\[\int_{\partial B_{r}}div(V(|x|)\nabla u)ds-\int_{\partial B}W(|x|)udx\geq 0,\]
and therefore,
\[V(r)\int_{\partial B_{r}}\Delta u dS -V_r(r)\int_{\partial B_{r}} \nabla u.x %%@
ds-W(r)\int_{\partial B_{r}}u(x)ds\geq 0.\]
It follows that 
\begin{equation}
V(r)\int_{\partial B_{r}}\Delta u dS -V_r(r)y'(r)-W(r)y(r)\geq 0, 
\end{equation}
and in view of (\ref{sup-eq}),  we see that $y$ satisfies the inequality
\begin{equation}\label{ode}
V(r)y''(r)+(\frac{(n-1)V(r)}{r}+V_r(r))y'(r)\leq -W(r)y(r), \ \ \ \ for\ \ \ 0<r<R, 
\end{equation}
that is it is a positive supersolution for $(B_{V,W})$. 

Standard results in ODE now allow us to conclude that $(B_{V,W})$ has actually a positive solution %%@
on $(0, R)$, and the proof of theorem \ref{main} is now complete. 
%\hfill $\square$.\\

%{\bf Proof of Theorem \ref{main-Rn}:} The proof is similar to that of Theorem \ref{main} so we %%@
%omit it. \hfill $\Box$

\subsection{Integral criteria for Bessel pairs}

In order to obtain criteria on $V$ and $W$ so that inequality $({\rm H}_{V, W})$ holds, we clearly %%@
need to investigate  whether   the ordinary differential equation $(B_{V,W})$ has positive %%@
solutions. For that, we rewrite $(B_{V,W})$ as
\[
(r^{n-1}V(r)y')'+r^{n-1}W(r)y=0, 
\]
and then by setting  $s=\frac{1}{r}$ and  $x(s)=y(r)$, we see that $y$ is a solution of %%@
$(B_{V,W})$ on an interval $(0,\delta)$ if and only if $x$ is a positive solution for the equation
\begin{equation}\label{s-ode}
\hbox{$(s^{-(n-3)}V(\frac{1}{s})x'(s))'+s^{-(n+1)}W(\frac{1}{s})x(s)=0$ \quad on \quad %%@
$(\frac{1}{\delta},\infty)$.}
\end{equation}
Now recall that a solution $x(s)$ of the equation (\ref{s-ode}) is said to be {\it oscillatory} if %%@
there exists a sequence $\{a_{n}\}^{\infty}_{n=1}$ such that $a_{n} \rightarrow +\infty$ and %%@
$x(a_{n})=0$. Otherwise we call the solution {\it non-oscillatory}. It follows from Sturm %%@
comparison theorem that all solutions of (\ref{s-ode}) are either all oscillatory or all %%@
non-oscillatory. Hence, the fact that $(V, W)$ is a Bessel pair or not 
is closely related to the oscillatory behavior of the equation (\ref{s-ode}). The following %%@
theorem is therefore a consequence of Theorem \ref{main}, combined with a relatively recent result %%@
of Sugie et al. in \cite{SKY} about the oscillatory behavior of the equation (\ref{s-ode}).

\begin{theorem}\label{main-cr} Let $V$ and $W$ be positive radial $C^1$-functions on 
$B_R\backslash \{0\}$, where $B_R$ is a ball centered at $0$ with radius $R$ in $\R^n$ ($n \geq 1$). Assume $\int^{R}_{0}\frac{1}{\tau^{n-1}V(\tau)}d\tau =+\infty$ and %%@
$\int_{0}^{R}r^{n-1}v(r)dr<\infty$. 
\begin{itemize}
\item Assume 
%there exists $0<\rho \leq R$ such that 
\begin{equation}\label{integral.1}
\limsup_{r\to 0}r^{2(n-1)}V(r)W(r)\big( \int^{R}_{r}\frac{1}{\tau^{n-1}V(\tau)}d\tau\big)^2< %%@
\frac{1}{4}
\end{equation}
 then  $(V, W)$ is a Bessel pair on $(0, \rho)$ for some $\rho>0$ and consequently,  inequality $({\rm H}_{V, W})$ %%@
holds for all $u \in C^{\infty}_{0}(B_{\rho})$, where $B_{\rho}$ is a ball of radius $\rho$. 
\item On the other hand, if
\begin{equation}\label{integral.2}
\liminf_{r\to 0}r^{2(n-1)}V(r)W(r)\big( \int^{R}_{r}\frac{1}{\tau^{n-1}V(\tau)}d\tau\big)^2> %%@
\frac{1}{4}
\end{equation}
then there is no interval $(0, \rho)$ on which $(V,W)$ is  a Bessel pair and consequently, there %%@
is no 
smooth domain $\Omega$ on which  inequality  
$({\rm H}_{V, W})$ holds.
\end{itemize}
\end{theorem}
A typical Bessel pair is $(|x|^{-\lambda}, |x|^{-\lambda -2})$ for $\lambda \leq n-2$. It is also easy to see by a simple change of variables in the corresponding ODEs that 
\begin{equation}
\hbox{$W$ is a Bessel potential if and only if  $
\left(|x|^{-\lambda}, |x|^{-\lambda}(|x|^{-2}+W(|x|)\right)$
is a Bessel pair.}
\end{equation}

More %%@
generally,   the above integral criterium allows to show the following.

\begin{theorem} \label{main.Bessel.pair} Let $V$ be an strictly positive $C^1$-function on $(0,R)$ such %%@
that for some $\lambda \in \R$ 
\begin{equation}
\hbox{$\frac{rV_r(r)}{V(r)}+\lambda \geq 0$ on $(0, R)$ and $\lim\limits_{r\to %%@
0}\frac{rV_r(r)}{V(r)}+\lambda =0$.}
\end{equation}
%\begin{enumerate}\item 
 If $\lambda \leq n-2$, then for any Bessel potential $W$ on $(0, R)$, and any $c\leq \beta %%@
(W; R)$, 
 the couple $(V, W_{\lambda, c})$ is a Bessel pair, where 
 \begin{equation}
 W_{\lambda,c}(r)=V(r)((\frac{n-\lambda-2}{2})^2r^{-2}+cW(r)).
 \end{equation}
Moreover, $\beta \big(V, W_{\lambda,c}; R\big)=1$ for all $c\leq \beta (W; R)$.
 
%\end{enumerate}
\end{theorem}
We need the following easy lemma. 
\begin{lemma}\label{strict-lemma} Assume the equation 
\[y''+\frac{a}{r}y'+V(r)y=0,\]
has a positive solution on $(0,R)$, where $a\geq 1$ and $V(r)> 0$. Then $y$ is strictly decreasing %%@
on $(0,R)$.
\end{lemma}
{\bf Proof:} First observe that $y$ can not have a local minimum, hence it is either increasing or %%@
decreasing on $(0,\delta)$, for $\delta$ sufficiently small. Assume $y$ is increasing. Under this %%@
assumption if $y'(a)=0$ for some $a>0$, then $y''(a)=0$ which contradicts the fact that $y$ is a %%@
positive solution of the above ODE. So we have $\frac{y''}{y'}\leq-\frac{a}{r},$ thus,
 \[y'\geq \frac{c}{r^a}.\]
 Therefore, $x(r) \rightarrow - \infty$ as $r \rightarrow 0$ which is a contradiction. Since, $y$ %%@
can not have a local minimum it should be strictly decreasing on $(0,R)$. \hfill $\Box$

{\bf Proof of Theorem \ref{main.Bessel.pair}:} Write $\frac{V_r(r)}{V(r)}=-\frac{\lambda}{r}+f(r)$
 where $f(r)\geq 0$ on $(0,R)$ and  $\lim\limits_{r\rightarrow 0}r f(r)=0$. In order to prove that 
$\left(V(r), V(r)((\frac{n-\lambda-2}{2})^2r^{-2}+cW(r))\right)$ is a Bessel pair, we need to 
show that the equation
\begin{equation}\label{f-ode.1}
y''+(\frac{n-\lambda-1}{r}+f(r))y'+((\frac{n-\lambda-2}{2})^2r^{-2}+cW(r))y(r)=0,
\end{equation}\label{f-ode}
has a positive solution on $(0,R)$. But first we note that the equation 
\[x''+(\frac{n-\lambda-1}{r})x'+((\frac{n-\lambda-2}{2})^2r^{-2}+cW(r))x(r)=0,\]
has a positive solution on $(0,R)$, whenever $c\leq \beta (W; R)$. 
Since now $f(r)\geq 0$ and since, by the proceeding lemma, $x'(r)\leq 0$, we get that  $x$ is a %%@
positive subsolution for the equation (\ref{f-ode.1}) on $(0,R)$, and thus it has a positive %%@
solution of $(0,R)$.  Note that this means that $\beta (V, W_{\lambda, c}; R) \geq 1$. 

For the reverse inequality, we shall use the criterium in Theorem \ref{main-cr}. Indeed apply %%@
criteria (\ref{integral.1}) to $V(r)$ and $W_1(r)=C\frac{V(r)}{r^2}$ to get 
\begin{eqnarray*}
\lim_{r \rightarrow 0}r^{2(n-1)}V(r)W_1(r)\big( %%@
\int^{R}_{r}\frac{1}{\tau^{n-1}V(\tau)}d\tau\big)^2&=&C\lim_{r \rightarrow %%@
0}r^{2(n-2)}V^{2}(r)\big( \int^{R}_{r}\frac{1}{\tau^{n-1}V(\tau)}d\tau\big)^2\\
&=&C\big(\lim_{r \rightarrow 0}r^{(n-2)}V(r)\int^{R}_{r}\frac{1}{\tau^{n-1}V(\tau)}d\tau\big)^2\\
&=&C\big(\lim_{r\rightarrow 0}\frac{\frac{1}{r^{n-1}V(r)}} %%@
{\frac{(n-2)r^{n-3}V(r)+r^{n-2}V_r(r)}{r^{2(n-2)}V^{2}(r)}}\big)^2\\
&=&C\big(\lim_{r \rightarrow 0}\frac{1}{(n-2)+r\frac{V_r(r)}{V(r)}}\big)^2\\
&=&\frac{C}{(n-\lambda-2)^2}.
\end{eqnarray*}
For $\big(V, CV(r^{-2} +cW)\big)$ to be a Bessel pair, it is necessary that  
$
\frac{C}{(n-\lambda-2)^2} 
%\leq C\lim_{r \rightarrow 0}r^{2(n-1)}V(r)W_{\lambda, c}(r)\big( %%@
%\int^{R}_{r}\frac{1}{\tau^{n-1}V(\tau)}d\tau\big)^2
\leq \frac{1}{4},
$ 
and the proof for the best constant is complete. \hfill $\Box$

With a similar argument one can also prove the following.
\begin{corollary} Let $V$ and $W$ be positive radial $C^1$-functions  on $B_R\backslash %%@
\{0\}$, where $B_R$ is a ball centered at zero with radius $R$ in $\R^n$ ($n \geq 1$).  Assume %%@
that 
\begin{equation}
\hbox{$\lim\limits_{r
\rightarrow 0}r\frac{V_r(r)}{V(r)}=-\lambda$ and $\lambda \leq n-2$.}
\end{equation} 
\begin{itemize}
\item If $\limsup\limits_{r\rightarrow 0}r^{2}\frac{W(r)}{V(r)}<(\frac{n-\lambda-2}{2})^2$, then %%@
$(V, W)$ is a Bessel pair on some interval $(0, \rho)$, and consequently there exists a ball $B_\rho  %%@
\subset R^n$ such that inequality $({\rm H}_{V, W})$ holds for all $u \in C^{\infty}_{0}(B_\rho)$.
\item On the other hand, if  $\liminf\limits_{r\rightarrow %%@
0}r^{2}\frac{W(r)}{V(r)}>(\frac{n-\lambda-2}{2})^2$, then there is no smooth domain $\Omega %%@
\subset R^n$ such that inequality $({\rm H}_{V, W})$ holds  on $\Omega$. 
\end{itemize}
\end{corollary}

\subsection{New weighted Hardy inequalities}

An immediate application of Theorem \ref{main.Bessel.pair} and Theorem \ref{main} is the following very general Hardy inequality. 

\begin{theorem} \label{super.hardy} Let $V(x)=V(|x|)$ be a strictly positive radial function on a smooth domain $\Omega$ containing %%@
$0$ such that  $R=\sup_{x \in \Omega} |x|$. Assume 
%Let $V$ be an strictly positive $C^1$-function on $(0,R)$ such %%@
that for some $\lambda \in \R$ 
\begin{equation}
\hbox{$\frac{rV_r(r)}{V(r)}+\lambda \geq 0$ on $(0, R)$ and $\lim\limits_{r\to %%@
0}\frac{rV_r(r)}{V(r)}+\lambda =0$.}
\end{equation}
\begin{enumerate}
\item  If $\lambda \leq n-2$, then  the following  
inequality holds  for any Bessel potential $W$ on $(0, R)$:  
\begin{equation}\label{v-hardy}
\hbox{$\int_{\Omega}V(x)|\nabla u|^{2}dx\geq  
(\frac{n-\lambda-2}{2})^2\int_{\Omega}\frac{V(x)}{|x|^2}u^{2}dx+\beta (W; R)\int_{\Omega} %%@
V(x)W(x)u^{2}dx$\quad  for all $u \in C^{\infty}_{0}(\Omega)$,}
\end{equation}
and both $(\frac{n-\lambda-2}{2})^2$ and $\beta (W; R)$ are the best constants. 

\item In particular, $\beta (V, r^{-2}V; R)=(\frac{n-\lambda-2}{2})^2$ is the best constant in the %%@
following inequality
\begin{equation}\label{super-hardy}
\hbox{$\int_{\Omega}V(x)|\nabla u|^{2}dx\geq  
(\frac{n-\lambda-2}{2})^2\int_{\Omega}\frac{V(x)}{|x|^2}u^{2}dx$ \quad  for all $u \in %%@
C^{\infty}_{0}(\Omega)$.}
\end{equation}
 
\end{enumerate}
\end{theorem}
Applied to $V_1(r)=r^{-m} W_{k, \rho}(r)$ and $V_2(r)=r^{-m} {\tilde W}_{k, \rho}(r)$ where
$ W_{k, \rho} %%@
(r)=\Sigma_{j=1}^k\frac{1}{r^{2}}\big(\prod^{j}_{i=1}log^{(i)}\frac{\rho}{r}\big)^{-2}$ and %%@
$\tilde W_{k; \rho} %%@
(r)=\Sigma_{j=1}^k\frac{1}{r^{2}}X^{2}_{1}(\frac{r}{\rho})X^{2}_{2}(\frac{r}{\rho}) \ldots  
X^{2}_{j-1}(\frac{r}{\rho})X^{2}_{j}(\frac{r}{\rho})$ are the iterated logs  introduced in the %%@
introduction, and noting that in both cases the corresponding $\lambda$ is equal to $2m+2$, we get the following new Hardy inequalities. 
 
\begin{corollary} Let $\Omega$ be a smooth bounded domain in $\R^n$ ($n \geq 1$) and $m\leq %%@
\frac{n-4}{2}$. Then the following inequalities hold.
\begin{eqnarray} 
\int_{\Omega}\frac{W_{k,\rho}(x)}{|x|^{2m}}|\nabla %%@
u|^2dx\geq(\frac{n-2m-4}{2})^2\int_{\Omega}\frac{W_{k,\rho}(x)}{|x|^{2m+2}}u^2dx\\
\int_{\Omega}\frac{\tilde W_{k,\rho}(x)}{|x|^{2m}}|\nabla %%@
u|^2dx\geq(\frac{n-2m-4}{2})^2\int_{\Omega}\frac{\tilde W_{k,\rho}(x)}{|x|^{2m+2}}u^2dx.
\end{eqnarray}
Moreover, the constant $(\frac{n-2m-4}{2})^2$ is the best constant in both %%@
inequalities.
\end{corollary}

\begin{remark}\rm
The two following theorems deal with Hardy-type  inequalities on the whole of $\R^n$.
%which hold for all $u \in C^{\infty}_{0}(\R^n)$.  
 Theorem \ref{main} already yields that inequality $(H_{V,W})$  holds for all $u \in C^{\infty}_{0}(\R^n)$ if and only if the ODE $(B_{V,W})$ has a positive solution on $(0,\infty)$. The latter equation is therefore 
 % Hence, if $(H_{V,W})$ (with $B=\R^n$)  holds for all $u \in C^{\infty}_{0}(\R^n)$, then it is necessary that the ordinary differential equation $(B_{V,W})$ to be 
 non-oscillatory, which will again be a very useful fact for computing  best constants, in view of the following criterium at infinity (Theorem 2.1 in \cite{SKY}) applied to the  equation 
 % in the next two theorems.  To be more precise, consider the equation  
\begin{equation} \label{remark-ode}
\left(a(r)y' \right)'+b(r)y(r)=0,
\end{equation}
where $a(r)$ and $b(r)$ are positive real valued functions.  Assuming that  $\int^{\infty}_{d} \frac{1}{a(\tau)}d\tau<\infty $ for some $d>0$, and that the following limit
\[L:=\lim_{r\rightarrow \infty }a(r)b(r)\left( 
\int_{r}^{\infty}\frac{1}{a(r)}dr \right)^2,\]
exists. Then  for the equation (\ref{remark-ode}) equation to be non-oscillatory, it is necessary that  $L\leq \frac{1}{4}$. 
%This will be one of our main tools to find the best constant in the following two theorems.
\end{remark}

\begin{theorem} \label{GM-IV} Let $a,b> 0$, and $\alpha, \beta, m$ be real numbers. 
\begin{itemize} 
\item If $\alpha\beta>0$, and $m\leq \frac{n-2}{2}$, %or $\alpha, \beta<0$, and $2m-\alpha \beta %%@< n-2$, 
then   for all $u \in C^{\infty}_{0}(\R^n)$
\begin{equation}\label{GM-V1}
\int_{\R^n}\frac{(a+b|x|^{\alpha})^{\beta}}{|x|^{2m}}|\nabla u|^2dx\geq %%@
(\frac{n-2m-2}{2})^2\int_{\R^n}\frac{(a+b|x|^{\alpha})^{\beta}}{|x|^{2m+2}}u^2dx, 
\end{equation}
and $(\frac{n-2m-2}{2})^2$ is the best constant in the inequality.
\item If $\alpha \beta<0$,  and $2m-\alpha \beta \leq n-2$, 
%or $\alpha>0$, $\beta<0$, and %%@$2m< n-2$, 
then   for all $u \in C^{\infty}_{0}(\R^n)$
\begin{equation}\label{GM-V2}
\int_{\R^n}\frac{(a+b|x|^{\alpha})^{\beta}}{|x|^{2m}}|\nabla u|^2dx\geq (\frac{n-2m+\alpha %%@
\beta-2}{2})^2\int_{\R^n}\frac{(a+b|x|^{\alpha})^{\beta}}{|x|^{2m+2}}u^2dx, 
\end{equation}
and $(\frac{n-2m+\alpha %%@
\beta-2}{2})^2$ is the best constant in the inequality.
\end{itemize}
\end{theorem}
{\bf Proof:} Letting $V(r)=\frac{(a+br^\alpha)^\beta}{r^{2m}}$, then
\[r\frac{V'(r)}{V(r)}=-2m+\frac{b\alpha \beta r^\alpha}{a+br^{\alpha}}=-2m+\alpha \beta -\frac{a %%@
\alpha \beta}{a+br^\alpha}.\]
Hence, in the case $\alpha, \beta> 0$ and $2m\leq n-2$, (\ref{GM-V1}) follows directly from %%@
Theorem \ref{super.hardy}. The same holds for (\ref{GM-V2}) since it also follows directly from Theorem \ref{super.hardy} in the case  where
$\alpha <0$, $\beta >0$ and $2m-\alpha \beta \leq n-2$. 

For the remaining two other cases,  we will use Theorem \ref{main}. Indeed, in this case the equation $(B_{V,W})$ becomes 
\begin{equation}\label{V-ODE1}
y''+(\frac{n-2m-1}{r}+\frac{b\alpha \beta r^{\alpha-1}}{a+br^{\alpha}})y'+\frac{c}{r^2}y=0,
\end{equation}
and the best constant in  inequalities (\ref{GM-V1}) and (\ref{GM-V2}) is the largest $c$ such %%@
that the above equation has a positive solution on $(0,+\infty)$. Note that by Lemma %%@
\ref{strict-lemma}, we have that  $y'<0$ on $(0,+\infty)$. Hence, if $\alpha <0$ and $\beta<0$, then the positive %%@
solution of the equation 
\[y''+\frac{n-2m-1}{r}y'+\frac{(\frac{n-2m-2}{2})^2}{r^2}y=0\]
 is a positive super-solution for (\ref{V-ODE1}) and therefore the latter ODE has a positive %%@
solution on $(0,+\infty)$, from which we conclude that (\ref{GM-V1}) holds. To prove now that $(\frac{n-2m-2}{2})^2$ is the best constant in (\ref{GM-V1}), we use the fact that if the equation (\ref{V-ODE1}) has a positive 
solution on $(0, +\infty)$, then the equation is necessarily non-oscillatory. By rewriting 
% the equation %%@
(\ref{V-ODE1}) as
\begin{equation} \label{V-ODE2}
\left(r^{n-2m-1}(a+br^{\alpha})^{\beta} y'\right)'+cr^{n-2m-3}(a+br^{\alpha})^{\beta}y=0, 
\end{equation}
and by noting that 
\[\int_{d}^{\infty}\frac{1}{r^{n-2m-1}(a+br^{\alpha})^{\beta}}<\infty,\]
and
\[\lim_{r\rightarrow \infty }cr^{2(n-2m-2)}(a+br^{\alpha})^{2\beta}\left( %%@
\int_{r}^{\infty}\frac{1}{r^{n-2m-1}(a+br^{\alpha})^{\beta}}dr \right)^2=\frac{c}{(n-2m-2)^2}, \]
we can use Theorem 2.1 in \cite{SKY} to conclude that for equation (\ref{V-ODE2}) to be non-oscillatory it is necessary that
\[\frac{c}{(n-2m-2)^2}\leq \frac{1}{4}.\]
Thus, $\frac{(n-2m-2)^2}{4}$ is the best constant in the inequality (\ref{GM-V1}). 

A very similar argument applies in the case where   $\alpha>0$, $\beta<0$, and $2m<n-2$, to obtain that  inequality (\ref{GM-V2}) holds for all $u \in C_{0}^{\infty}(\R^n)$ and that  $(\frac{n-2m+\alpha \beta-2}{2})^2$  
is indeed the best constant. 
%The proof is complete.
 \hfill $\Box$\\

Note that the above two inequalities can be improved on smooth bounded domains by using Theorem \ref{super.hardy}.  

We shall now extend the recent results of Blanchet-Bonforte-Dolbeault-Grillo-Vasquez \cite{BBDGV} and address some of their questions regarding best constants.  

\begin{theorem}\label{GM-V} Let $a,b> 0$, and $\alpha, \beta$ be real numbers. 
\begin{itemize}
\item If $\alpha \beta<0$ and $-\alpha \beta \leq n-2$, %or $\alpha>0$, $\beta<0$, and $n\geq %%@2$, 
then   for all $u \in C^{\infty}_{0}(\R^n)$
\begin{equation}\label{GM-V3}
\int_{\R^n}(a+b|x|^{\alpha})^{\beta}|\nabla u|^2dx\geq b^{\frac{2}{\alpha}}(\frac{n-\alpha %%@
\beta-2}{2})^2\int_{\R^n}(a+b|x|^{\alpha})^{\beta-\frac{2}{\alpha}}u^2dx,
\end{equation}
and $b^{\frac{2}{\alpha}}(\frac{n-\alpha 
\beta-2}{2})^2$ is the best constant in the inequality.
\item If $\alpha \beta >0$ and %$-\alpha \beta \leq n-2 $, or $\alpha>0$, $\beta>0$, and %%@
$n\geq 2$, then there exists a constant $C>0$ such that   for all $u %%@
\in C^{\infty}_{0}(\R^n)$
\begin{equation}\label{GM-V4}
\int_{\R^n}(a+b|x|^{\alpha})^{\beta}|\nabla u|^2dx\geq C %%@
\int_{\R^n}(a+b|x|^{\alpha})^{\beta-\frac{2}{\alpha}}u^2dx.
\end{equation}
Moreover, $b^{\frac{2}{\alpha}}(\frac{n-2}{2})^2\leq C\leq b^{\frac{2}{\alpha}}(\frac{n+\alpha %%@
\beta-2}{2})^2$.
\end{itemize}
\end{theorem}
{\bf Proof:}  Letting $V(r)=(a+br^{\alpha})^{\beta}$, then we have
\[
r\frac{V'(r)}{V(r)}=\frac{b\alpha \beta r^{\alpha}}{a+br^{\alpha}}=\alpha \beta-\frac{a \alpha %%@
\beta}{a+br^{\alpha}}.
\]
Inequality (\ref{GM-V3}) and its best constant in the case when $\alpha<0$ and $\beta>0$, then follow immediately from Theorem \ref{super.hardy} with $\lambda=-\alpha \beta$.
The proof of the remaining cases will use Theorem \ref{main} as well as the integral criteria for %%@
the oscillatory behavior of solutions for ODEs of the form ($B_{V,W}$). 
%Indeed, letting $V(r)=(a+br^{\alpha})^{\beta}$, then we have
%\[r\frac{V'(r)}{V(r)}=\frac{b\alpha \beta r^{\alpha}}{a+br^{\alpha}}=\alpha \beta-\frac{a \alpha %%@
%\beta}{a+br^{\alpha}}.\]

Assuming still that $\alpha \beta <0$, then with an argument similar to that of Theorem \ref{GM-IV} above, one can show %%@
that the positive solution of the equation $y''+(\frac{n+\alpha \beta-1}{r})y'+\frac{(n+\alpha %%@
\beta -2)^2}{4r^2}y=0$ on $(0,+\infty)$ is a positive supersolution for the equation
\[y''+(\frac{n-1}{r}+\frac{V'(r)}{V(r)})y'+\frac{b^{\frac{2}{\alpha}}(n+\alpha \beta %%@
-2)^2}{4(a+br^{\alpha})^{\frac{2}{\alpha}}}y=0.\]
Theorem \ref{main} then yields that the inequality (\ref{GM-V3}) holds for all $u \in %%@
C_{0}^{\infty}(\R^n)$.
To prove now that $b^{\frac{2}{\alpha}}(\frac{n+\alpha \beta -2}{2})^2$ is the best constant in %%@
(\ref{GM-V3}) it is enough to show that if the following equation
\begin{equation}\label {V-ODE3}
\left(r^{n-1}(a+br^{\alpha})^{\beta}y'\right)'+cr^{n-1}(a+br^{\alpha})^{\beta-\frac{2}{\alpha}}y=0
\end{equation}
has a positive solution on $(0,+\infty)$, then $c\leq b^{\frac{2}{\alpha}}(\frac{n+\alpha \beta %%@
-2}{2})^2$. 
%Here we consider two different case. If $\alpha<0$ and $\beta>0$, then we have
%\[ \lim_{r \rightarrow %%@
%0}cr^{2(n-1)}(a+br^{\alpha})^{2\beta-\frac{2}{\alpha}}\left(\int_{r}^{a}\frac{1}{r^{n-1}(a+br^{
%\alpha})^{\beta}}dr\right)^2=\frac{c}{b^{\frac{2}{\alpha}}(n+\alpha \beta-2)^2}.\]
%Hence, it follows from Theorem \ref{main-cr} that $c\leq \frac{(n+\alpha \beta-2)^2}{4}$.
 If  now $\alpha>0$ and $\beta<0$, then we have 
\[\lim_{r\rightarrow \infty }cr^{2(n-1)}(a+br^{\alpha})^{2\beta-\frac{2}{\alpha}}\left( %%@
\int_{r}^{\infty}\frac{1}{r^{n-1}(a+br^{\alpha})^{\beta}}dr %%@
\right)^2=\frac{c}{b^{\frac{2}{\alpha}}(n+\alpha \beta-2)^2}.\] 
Hence, by Theorem 2.1 in \cite{SKY} again, the non-oscillatory aspect of the equation holds for $c\leq \frac{b^{\frac{2}{\alpha}}(n+\alpha %%@
\beta-2)^2}{4}$ which completes the proof of the first part. \\
A similar argument applies in the case where $\alpha \beta>0$ to prove that  (\ref{GM-V4}) holds for all $u \in C_{0}^{\infty}(\R^n)$ and %%@
$b^{\frac{2}{\alpha}}(\frac{n-2}{2})^2\leq C\leq b^{\frac{2}{\alpha}}(\frac{n+\alpha %%@
\beta-2}{2})^2$.  
The  best constants are estimated by carefully studying  the existence of positive solutions for the ODE (\ref{V-ODE3}).

\begin{remark}\rm Recently, Blanchet et al. in \cite{BBDGV} studied a special case of  inequality %%@
(\ref{GM-V3}) ($a=b=1$, and $\alpha=2$) under the additional condition:
\begin{equation}\label{extra-cond}
\int_{\R^n} (1+|x|^2)^{\beta-1}u(x)dx=0, \ \ for \ \ \beta<\frac{n-2}{2}.
\end{equation}
Note that  we do not assume (\ref{extra-cond}) in Theorem \ref{GM-V}, and that we have found the best constants for %%@
$\beta\leq 0$, a case that was left open  in \cite{BBDGV}.
% optimality of the constants has been left open for all %%@$\beta\leq 0$. 
\end{remark}

 \subsection{Improved Hardy and Caffarelli-Kohn-Nirenberg Inequalities}

%It is clear that Theorem \ref{super.hardy} is a major extension of the classical Hardy's %%@
%inequality. Indeed, applied with $V(r)\equiv 1$, we get that for any  domain $\Omega$ in $R^{n}$, %%@
%$n \geq 3$,  we have
%\begin{equation}\label{cl-hardy}
%\hbox{$\int_{\Omega}|\nabla u |^{2}dx \geq ( \frac{n-2}{2})^{2} %%@
%\int_{\Omega}\frac{|u|^2}{|x|^{2}}dx$ \quad  for all $u \in H^{1}_{0}(\Omega)$.}
%\end{equation}
 In \cite{CKN} Caffarelli-Kohn-Nirenberg established a  set inequalities of the %%@
following form:
\begin{equation}\label{CKN}
\hbox{
$\big(\int_{R^n}|x|^{-bp}|u|^{p}dx\big)^{\frac{2}{p}}\leq C_{a,b}\int_{R^n}|x|^{-2a}|\nabla %%@
u|^2dx$ for all $u \in C^{\infty}_{0}(R^n)$,}
\end{equation}
where for $n\geq 3$,
\begin{equation}\hbox{
$-\infty<a<\frac{n-2}{2}$, $a\leq b \leq a+1,$ and $p=\frac{2n}{n-2+2(b-a)}$.
}
\end{equation}
For the cases $n=2$ and $n=1$ the conditions are slightly different. For $n=2$
\begin{equation}\hbox{
$-\infty<a<0$, $a< b \leq a+1,$ and $p=\frac{2}{b-a}$,
}
\end{equation}
and for $n=1$ 
\begin{equation}\hbox{
$-\infty<a<-\frac{1}{2}$, $a+\frac{1}{2}< b \leq a+1,$ and $p=\frac{2}{-1+2(b-a)}$.
}
\end{equation}
Let $D_{a}^{1,2}$ be the completion of $C^{\infty}_{0}(R^n)$ for the inner product 
$(u,v)=\int_{R^n}|x|^{-2a}\nabla u. \nabla vdx$ and 
let 
\begin{equation}
S(a,b)=\inf_{u \in D_{a}^{1,2}\backslash \{0\}} \frac{\int_{R^n}|x|^{-2a}|\nabla %%@
u|^2dx}{(\int_{R^n}|x|^{-bp}|u|^{p}dx\big)^{2/p}}
\end{equation}
denote the best embedding constant.
We are concerned here with the ``Hardy critical" case of the above inequalities, that is when %%@
$b=a+1$.  In this direction, Catrina and Wang \cite{CW}   showed that for $n\geq 3$    we have %%@
$S(a,a+1)=(\frac{n-2a-2}{2})^2$ and that $S(a,a+1)$ is not achieved while $S(a,b)$ is always %%@
achieved for $a<b<a+1$. For the case $n=2$ they also showed that   $S(a,a+1)=a^{2}$, and that %%@
$S(a,a+1)$ is not achieved, while for $a<b<a+1$,  $S(a,b)$ is again achieved. For $n=1$,  %%@
$S(a,a+1)=(\frac{1+2a}{2})^{2}$ is also not achieved. 

In this section we give a necessary and sufficient condition for improvement of (\ref{CKN}) with %%@
$b=a+1$ and $n\geq 1$. Our results cover also the critical case when $a=\frac{n-2}{2}$ which is %%@
not allowed by the methods of  \cite{CKN}. 

\begin{theorem} \label{main-CKN} Let $W$ be a positive radial function on the ball $B$ in $\R^n$ %%@
($n \geq 1$) with radius $R$ and centered at zero.
Assume $a\leq \frac{n-2}{2}$. The  following two statements are then equivalent:

\begin{enumerate}
\item $W$ is a Bessel potential on $(0, R)$. 
 
\item There exists $c>0$ such that the following inequality holds for all $u \in %%@
C_{0}^{\infty}(B)$
\begin{equation*}
\hbox{$({\rm H}_{a, cW})$ \quad \quad \quad \quad \quad  $\int_{B}|x|^{-2a}|\nabla u(x) |^{2}dx %%@
\geq (\frac{n-2a-2}{2})^2\int_{B}|x|^{-2a-2}u^2 dx+c\int_{B} |x|^{-2a}W(x)u^2dx,$\quad \quad \quad %%@
\quad \quad \quad \quad  \quad \quad \quad \quad }
\end{equation*} 
\end{enumerate}
 Moreover, $(\frac{n-2a-2}{2})^2$ is the best constant and $\beta (W; R)=\sup \{c; (H_{a,cW}) %%@
holds\}$, where $\beta (W; R)$ is the weight of the Bessel potential $W$ on $(0, R)$. 

On the other hand, there is no strictly positive $W \in C^{1}(0,\infty)$, 
%and $W(r)>0$ for all $r \in (0, \infty)$. Then %%@ there is no $c>0$ 
such that the following inequality holds for all $u \in C_{0}^{\infty}(\R^n)$, 
\begin{equation}\label{no-improve}
\int_{\R^n}|x|^{-2a}|\nabla u(x) |^{2}dx \geq (\frac{n-2a-2}{2})^2\int_{\R^n}|x|^{-2a-2}u^2 %%@
dx+c\int_{\R^n} W(|x|)u^2dx.
\end{equation}

 \end{theorem}
 {\bf Proof:} It suffices to use Theorems  \ref{main} and  \ref{super.hardy} with $V(r)=r^{-2a}$ %%@
to get that $W$ is a Bessel function if and only if the pair $\big(r^{-2a}, W_{a,c}(r)\big)$  is a %%@
Bessel pair on $(0, R)$ for some $c>0$, where 
 \[W_{a,c}(r)=(\frac{n-2a-2}{2})^2r^{-2-2a}+cr^{-2a}W(r).\]
 For the last part, assume that (\ref{no-improve}) holds for some $W$. Then it follows from %%@
Theorem \ref{main-CKN} that for $V=cr^{2a}W(r)$ the equation  $y''(r)+\frac{1}{r}y' +v(r)y=0$ has %%@
a positive solution on $(0,\infty)$. From Lemma \ref{strict-lemma} we know that $y$ is strictly %%@
decreasing on $(0,+\infty)$. Hence, $\frac{y''(r)}{y'(r)}\geq -\frac{1}{r}$ which yields $y'(r)\leq %%@
\frac{b}{r}$, for some $b>0$. Thus $y(r) \rightarrow -\infty$ as $r\rightarrow +\infty$. This is a %%@
contradiction and the proof is complete. \hfill $\Box$

\begin{remark} \rm Theorem \ref{main-CKN} characterizes the best constant only when $\Omega$ is a ball, while for general domain $\Omega$,  it just gives a lower and upper bounds for the best constant corresponding to a given Bessel potential $W$. It is indeed clear that 
\[C_{B_{R}}(W) \leq C_{\Omega}(W)\leq C_{B_{\rho}}(W),\]
where $B_{R}$ is the smallest ball containing $\Omega$ and $B_{\rho}$ is the largest %%@
ball contained in it. If now $W$ is a Bessel potential such that $\beta (W, R)$ is independent of $R$, then clearly $\beta (W, R)$ is also the best constant  in inequality $(H_{a,cW})$ for any smooth bounded domain.
%$C_{B_{R}}(V)=C_{B_{\rho}}(V)$ for some potential %%@
%$W$ then we will be able to find the best constant $C_{\Omega}(V)$ in inequality $(H_{a,cW})$ for %%@
%the non-radial domain $\Omega$. This is why it is easy to find the best constant for the %%@
This is clearly the case for the  potentials $W_{k, \rho}$ and 
%=\frac{1}{4|x|^2}\Sigma^{k}_{j=1} (\prod^{j}_{i=1}\log^{(i)}\frac{\rho}{|x|})^{-2}$ and  \[
$\tilde W_{k, \rho}$
%=\frac{1}{4|x|^2} \Sigma^{k}_{j=1} X^{2}_{1}(\frac{|x|}{D})X^{2}_ {2}(\frac{|x|}{D}) \ldots %%@
%X^{2}_{j-1}(\frac{|x|}{D})X^{2}_{j}(\frac{|x|}{D}),
%\]
where $\beta (W, R)=\frac{1}{4}$ for all $R$, while for $W\equiv 1$ the best constant is still not known for general domains even for the simplest case %%@
$a=0$. 
\end{remark}
Using the integral criteria for Bessel potentials, we can also  deduce immediately the following. 

\begin{corollary} Let $\Omega$ be a bounded smooth domain in $R^n$ with $n\geq 1$, and let $W$ be %%@
a non-negative  function in $C^{1}(0, R=:\sup_{x \in \partial \Omega}|x|]$ and $a\leq %%@
\frac{n-2}{2}$.   
\begin{enumerate}
\item If 
$\hbox{$\liminf\limits_{r\rightarrow 0} \ln(r)\int^{r}_{0} sW(s)ds>-\infty$, } $
then there exists $\alpha:=\alpha(\Omega)>0$ such that an improved Hardy inequality $({\rm H}_{a, %%@
W_{\alpha}})$ holds for the scaled potential $W_\alpha(x):=\alpha^2W(\alpha |x|)$.  
\item  If
$\hbox{$\lim\limits_{r\rightarrow 0} \ln(r)\int^{r}_{0} sW(s)ds=-\infty$, } $
then there are no $\alpha, c>0$,  for which $({\rm H}_{a,W_{\alpha,c}})$ holds with %%@
$W_{\alpha,c}=cW(\alpha |x|)$.
\end{enumerate}
\end{corollary}

By applying the above to various examples of Bessel potentials, we can now deduce several old and %%@
new inequalities.  The first is an extension of a result established by  Brezis and V\'{a}zquez %%@
\cite{BV} in the case where $a=0$, and $b=0$. 

\begin{corollary} \label{CKN.BV} Let $\Omega$ be a bounded smooth domain in $R^n$ with $n\geq 1$ %%@
and $a\leq \frac{n-2}{2}$. Then, for any $b<2a+2$ there exists $c>0$ such that for all $u \in %%@
C_{0}^{\infty}(\Omega)$
\begin{equation}\label{beta}
\hbox{$\int_{\Omega}|x|^{-2a}|\nabla u |^{2}dx \geq %%@
(\frac{n-2a-2}{2})^2\int_{\Omega}|x|^{-2a-2}u^2 dx+c\int_{\Omega} |x|^{-b}u^2dx.$}
\end{equation} 
Moreover, when $\Omega$ is a ball $B$ of radius $R$ the best constant $c$ for which (\ref{beta})  %%@
holds is equal to the weight $\beta (r^{2a-b}; R)$ of the Bessel potential $W(r)=r^{2a-b}$ on $(0, %%@
R]$.\\
 In particular, 
\begin{equation}
\hbox{ \quad \quad \quad \quad \quad \quad  $\int_{B}|x|^{-2a}|\nabla u|^{2}dx \geq %%@
(\frac{n-2a-2}{2})^2\int_{B}|x|^{-2a-2}u^2 dx+\lambda_{B}\int_{B} %%@
|x|^{-2a}u^2dx,$\quad \quad \quad \quad \quad \quad \quad  \quad \quad \quad \quad \quad \quad}
\end{equation}
where the best constant $\lambda_{B}$ is equal to $z_{0}\omega^{2/n}_{n}|\Omega|^{-2/n}$,  %%@
where $\omega_{n}$ and $|\Omega|$ denote the volume of the unit ball and $\Omega$ respectively, %%@
and $z_{0}=2.4048...$ is the first zero of the Bessel function $J_{0}(z)$.
\end{corollary} 

{\bf Proof:}  It suffices to apply  Theorem \ref{main-CKN} with  the function $W(r)=r^{b+2a}$ %%@
which is a Bessel potential whenever $b >-2a-2$ since then 
  ${\rm liminf}_{r\rightarrow  0} \ln(r)\int^{r}_{0} s^{2a+1}W(s)ds>-\infty$.
In the case where $b=-2a$ and therefore $W\equiv 1$,  we use the fact that $\beta (1;  %%@
R)=\frac{z^2_0}{R^2}$ (established in the appendix) to deduce that  the best constant is then %%@
equal to $z_{0}\omega^{2/n}_{n}|\Omega|^{-2/n}$.
\hfill $\square$   \\

The following corollary is an extension of a recent result by Adimurthi et all \cite{ACR} %%@
established  in the case where $a=0$, and of another result by Wang and Willem in \cite{WW} %%@
(Theorem 2) in the case $k=1$. We also provide  here the value of the best constant.
\begin{corollary}\label{CKN.A} Let $B$ be a bounded smooth domain in $R^n$ with $n\geq 1$ and %%@
$a\leq \frac{n-2}{2}$. Then for every integer $k$,  and $\rho=(\sup_{x \in \Omega}|x|)( %%@
e^{e^{e^{.^{.^{e((k-1)-times)}}}}} )$, we have  for any $u \in H^{1}_{0}(\Omega)$, 
\begin{equation}\label{ar-hardy}
\hbox{$\int_{\Omega}|x|^{-2a}|\nabla u|^{2}dx \geq %%@
(\frac{n-2a-2}{2})^2\int_{\Omega}\frac{u^2}{|x|^{2a+2}} %%@
dx+\frac{1}{4}\sum^{k}_{j=1}\int_{\Omega}\frac{|u |^{2}}{|x|^{2a+2}}\big( %%@
\prod^{j}_{i=1}log^{(i)}\frac{\rho}{|x|}\big)^{-2}dx$.}
\end{equation}
 Moreover,  $\frac{1}{4}$ is the best constant  which is not attained in $H_{0}^{1}(\Omega)$.
\end{corollary}

{\bf Proof:}  As seen in the appendix, $W_{k, \rho}(r)=\sum^{k}_{j=1}\frac{1}{r^2}\big( %%@
\prod^{j}_{i=1}log^{(i)}\frac{\rho}{|x|}\big)^{-2}dx$ is a Bessel potential 
on $(0, R)$ where $R=\sup_{x \in \Omega}|x|$, and $\beta (W_{k, \rho}; R)=\frac{1}{4}$. 
\hfill $\square$\\

The very same reasoning leads to the following extension of a result established by Filippas and %%@
Tertikas    \cite{FT} in the case where $a=0$. 
\begin{corollary} \label{CKN.FT} Let $\Omega$ be a bounded smooth domain in $R^n$ with $n\geq 1$ %%@
and $a\leq \frac{n-2}{2}$. Then for every integer $k$,  and any $D\geq \sup_{x \in \Omega}|x|$, we %%@
have  for $u \in H^{1}_{0}(\Omega)$,
\begin{equation}\label{ar-hardy}
\hbox{$\int_{\Omega}\frac{|\nabla u|^{2}}{|x|^{2a}}dx \geq %%@
(\frac{n-2a-2}{2})^2\int_{\Omega}\frac{u^2}{|x|^{2a+2}} %%@
dx+\frac{1}{4}\sum^{\infty}_{i=1}\int_{\Omega}\frac{1}{|x|^{2a+2}}X
^{2}_{1}(\frac{|x|}{D})X^{2}_{2}(\frac{|x|}{D})...X^{2}_{i}(\frac{|x|}{D})|u|^{2}dx, $
}
\end{equation}
and  $\frac{1}{4}$ is the best constant  which is not attained in $H_{0}^{1}(\Omega)$.
 
\end{corollary}
%The following theorem shows that we can not improve inequality 
%\[\int_{\R^n}|x|^{-2a}|\nabla u(x) |^{2}dx \geq (\frac{n-2a-2}{2})^2\int_{\R^n}|x|^{-2a-2}u^2 %%@
%dx,\]
%with any strictly positive potential. 
%\begin{theorem} Assume $W \in C^{1}(0,\infty)$, and $W(r)>0$ for all $r \in (0, \infty)$. Then %%@
%there is no $c>0$ such that the following inequality holds for all $u \in C_{0}^{\infty}(\R^n)$.
%\begin{equation}\label{no-improve}
%\int_{\R^n}|x|^{-2a}|\nabla u(x) |^{2}dx \geq (\frac{n-2a-2}{2})^2\int_{\R^n}|x|^{-2a-2}u^2 %%@
%dx+c\int_{\R^n} W(|x|)u^2dx.
%\end{equation}
%\end{theorem}
%{\bf Proof:} Assume that (\ref{no-improve}) holds for some $W$ and $c>0$. Then it follows from %%@
%Theorem \ref{main-CKN} that for $V=cr^{2a}W(r)$ the equation  $y''(r)+\frac{1}{r}y' +v(r)y=0$ has %%@
%a positive solution on $(0,\infty)$. From Lemma \ref{strict-lemma} we know that $y$ is strictly %%@
%decreasing on $(0,\infty)$. Hence, $\frac{y''(r)}{y'(r)}\geq -\frac{1}{r}$ which yields $y'(r)\leq %%@
%\frac{b}{r}$, for some $b>0$. Thus $y(r) \rightarrow -\infty$, as $r\rightarrow \infty$. This is a %%@
%contradiction and the proof is complete. \hfill $\Box$

% \subsection{Hardy Inequalities in dimension two}
The classical Hardy inequality is valid for dimensions $n\geq 3$. We now present optimal Hardy %%@
type inequalities for dimension two in bounded domains, as well as the corresponding best constants.
\begin{theorem}\label{2dim-hardy} Let $\Omega$ be a smooth domain in $R^2$ and $0 \in \Omega$. %%@
Then we have the following inequalities.
\begin{itemize}
\item Let $D\geq \sup_{x \in \Omega}|x|$, then for all $u \in H^{1}_{0}(\Omega),$
\begin{equation}\hbox{$\int_{\Omega}|\nabla u |^{2}dx \geq %%@
\frac{1}{4}\sum^{\infty}_{i=1}\int_{\Omega}\frac{1}{|x|^{2}}X
^{2}_{1}(\frac{|x|}{D})X^{2}_{2}(\frac{|x|}{D})...X^{2}_{i}(\frac{|x|}{D})|u|^{2}dx$}
\end{equation}
and $\frac{1}{4}$ is the best constant.

\item Let $\rho=(\sup_{x \in \Omega}|x|)( e^{e^{e^{.^{.^{e((k-1)-times)}}}}} )$, then for all $u %%@
\in H^{1}_{0}(\Omega)$
\begin{equation}
\hbox{$\int_{\Omega}|\nabla u |^{2}dx \geq \frac{1}{4}\sum^{k}_{j=1}\int_{\Omega}\frac{|u %%@
|^{2}}{|x|^{2}}\big( \prod^{j}_{i=1}log^{(i)}\frac{\rho}{|x|}\big)^{-2}dx$,}
\end{equation} 
and $\frac{1}{4}$ is the best constant  for all $k\geq 1$.
\item If $\alpha<2$, then there exists $c>0$ such that for all $u \in H^{1}_{0}(\Omega),$
\begin{equation}
\hbox{$\int_{\Omega}|\nabla u |^{2}dx \geq c\int_{\Omega}\frac{u^2}{|x|^{\alpha}}\, dx$,}
\end{equation}
and the best constant  is larger or equal to $\beta (r^\alpha; \sup\limits_{x \in \Omega}|x|)$.

\end{itemize}

\end{theorem}
 
An immediate application of Theorem \ref{main} coupled with H\"older's inequality gives the %%@
following duality statement, which should  be compared to  inequalities  dual to those of %%@
Sobolev's, recently obtained via the theory of mass transport  \cite{AGK, CNV}.

\begin{corollary} \label{dual} Suppose that $\Omega$ is a smooth bounded domain containing $0$ in $R^{n}$ ($n\geq1$) with $R:=\sup_{x\in \Omega}|x|$. Then,  for any $a\leq \frac{n-2}{2}$ and $0<p\leq 2$, we have the following dual inequalities:
\begin{eqnarray*}
\inf \left\{\int_{\Omega}|x|^{-2a}|\nabla u |^{2}dx - ( \frac{n-2a-2}{2})^{2} %%@
\int_{\Omega}|x|^{-2a-2}|u|^{2}dx;\,  u \in C_{0}^{\infty}(\Omega), ||u||_p=1\right\}\\
 \geq \sup \left\{\big(\int_\Omega (\frac{|x|^{-2a}}{W(x)})^{\frac{p}{p-2}}\, dx %%@
\big)^{\frac{2-p}{p}}
; \, W\in  {\cal B}(0,R)\right\}.
\end{eqnarray*}

\end{corollary}
\section{General Hardy-Rellich inequalities}

Let $0 \in \Omega \subset R^n$ be a smooth domain, and denote
\[C^{k}_{0,r}(\Omega)=\{v \in C^{k}_{0}(\Omega): \mbox{v is radial and supp v }\subset \Omega\},\]
\[H^{m}_{0,r}(\Omega)=\{ u \in H^{m}_{0}(\Omega): \mbox{u is radial}\}.\]
We start by considering  a general inequality for radial functions.

\begin{theorem} \label{mainrad.hr} Let $V$ and $W$ be positive radial $C^1$-functions on a ball %%@
$B\backslash \{0\}$, where $B$ is a ball with radius $R$ in $\R^n$ ($n \geq 1$) and centered at %%@
zero. Assume $\int^{R}_{0}\frac{1}{r^{n-1}V(r)}dr=\infty$ and $\lim_{r \rightarrow %%@
0}r^{\alpha}V(r)=0$ for some $\alpha< n-2$. Then the following statements are equivalent: 
\begin{enumerate}

\item $(V, W)$ is a Bessel pair on $(0, R)$. 
\item There exists $c>0$ such that the  following inequality holds for  all radial functions $u %%@
\in  
C^{\infty}_{0,r}(B)$
\begin{equation*}\label{gen-hardy}
\hbox{   $({\rm HR}_{V,cW})$ \quad \quad \quad \quad  $\int_{B}V(x)|\Delta u |^{2}dx \geq  %%@
c\int_{B} W(x)|\nabla u|^{2}dx+(n-1)\int_{B}(\frac{V(x)}{|x|^2}-\frac{V_r(|x|)}{|x|})|\nabla %%@
u|^2dx.$ \quad \quad \quad \quad \quad \quad \quad \quad  }
\end{equation*}

\end{enumerate} 
Moreover,  the best constant is given by 
\begin{equation}
\hbox{$\beta (V, W; R)=\sup\big\{c; \, \,  ({\rm HR}_{V, cW})$ holds for radial %%@
functions$\big\}$.} 
\end{equation}

\end{theorem} 
{\bf Proof:}  Assume $u \in C^{\infty}_{0,r}(B)$ and observe that
\[\int_{B}V(x)|\Delta u %%@
|^{2}dx=n\omega_{n}\{\int^{R}_{0}V(r)u_{rr}^{2}r^{n-1}dr+(n-1)^2\int^{R}_{0}V(r)\frac{u^{2}_{r}}{r
^{2
}}r^{n-1}dr
+2(n-1)\int^{R}_{0}V(r)uu_rr^{n-2}dr\}.\]
Setting $\nu=u_{r}$, we then have
\[\int_{B}V(x)|\Delta u |^{2}dx=\int_{B}V(x)|\nabla \nu |^{2}dx+(n-1) %%@
\int_{B}(\frac{V(|x|)}{|x|^2}-\frac{V_r(|x|)}{|x|})|\nu|^{2}dx. \]
Thus, $({\rm HR}_{V,W})$ for radial functions is equivalent to
\[\int_{B}V(x)|\nabla \nu |^{2}dx\geq \int_{B}W(x)\nu^2 dx.\]
Letting $x(r)=\nu(x)$ where $|x|=r$, we then have
\begin{equation}\label{1-dim}
\int^{R}_{0}V(r)(x'(r))^2r^{n-1}dr \geq \int^{R}_{0}W(r)x^{2}(r)r^{n-1}dr.
\end{equation}
It therefore follows from Theorem \ref{main} that 1) and 2) are  equivalent. \hfill $\Box$\\

By applying the above theorem  to the Bessel pair 
\[
\hbox{$V(x)=|x|^{-2m}$ \quad and \quad $W_{m}(x)= 
V(x)\big[(\frac{n-2m-2}{2})^2|x|^{-2}+W(x)\big]$}
\]
 where $W$ is a Bessel potential, and by using Theorem \ref{super.hardy},   we get the following %%@
result  in the case of radial functions.  

\begin{corollary} \label{radial} Suppose $n\geq 1$ and $m<\frac{n-2}{2}$. Let $B_{R}\subset \R^n$ %%@
be a ball of radius $R>0$ and centered at zero. Let $W$ be a Bessel potential on $(0, R)$. Then we %%@
have  for all $u \in C^{\infty}_{0,r}(B_{R})$
\begin{equation}
\int_{B_{R}}\frac{|\Delta u|^2}{|x|^{2m}}\geq (\frac{n+2m}{2})^{2}\int_{B_{R}}\frac{|\nabla %%@
u|^2}{|x|^{2m+2}}dx+\beta (W; R)\int_{B_{R}}W(x)\frac{|\nabla u|^2}{|x|^{2m}}dx.
\end{equation}
Moreover,   $(\frac{n+2m}{2})^{2}$ and $\beta (W; R)$ are the best constants.
% in the following %%@
%sense:
%\begin{equation}
%(\frac{n+2m}{2})^{2}=\inf\left\{\frac{\int_{B_{R}}\frac{|\Delta %%@
%u|^2}{|x|^{2m}}dx-\int_{B_{R}}W(x)\frac{|\nabla u|^2}{|x|^{2m}}dx}{\int_{B_{R}}\frac{|\nabla %%@
%u|^2}{|x|^{2m+2}}dx};\, u \in C^{\infty}_{0,r}(B_{R})\setminus \{0\}\right\}
%\end{equation}
%and
%\begin{equation}
%\beta (W; R)=\inf\left\{ \frac{\int_{B_{R}}\frac{|\Delta u|^2}{|x|^{2m}}- %%@
%(\frac{n+2m}{2})^{2}\int_{B_{R}}\frac{|\nabla u|^2}{|x|^{2m+2}}dx}{\int_{B_{R}}W(x)\frac{|\nabla %%@
%u|^2}{|x|^{2m}}dx};\, u \in C^{\infty}_{0,r}(B_{R})\setminus \{0\}\right\}.
%\end{equation}

\end{corollary}

\subsection{The non-radial case}

The decomposition of a function into its spherical harmonics will be one of our tools to prove the %%@
corresponding  result in the non-radial case. This idea has also been used in \cite{TZ}. Any %%@
function $u \in C^{\infty}_{0}(\Omega)$ could be extended by zero outside $\Omega$, and could %%@
therefore be considered as a function in $C^{\infty}_{0}(R^n)$. By decomposing $u$ into spherical %%@
harmonics we get
\[
\hbox{$u=\Sigma^{\infty}_{k=0}u_{k}$ 
where $u_{k}=f_{k}(|x|)\varphi_{k}(x)$}
\]
 and $(\varphi_k(x))_k$ are the orthonormal eigenfunctions of the Laplace-Beltrami operator  %%@
with corresponding eigenvalues $c_{k}=k(n+k-2)$, $k\geq 0$. The functions $f_{k}$ belong to %%@
$C_{0}^{\infty}(\Omega)$ and satisfy $f_{k}(r)=O(r^k)$ and $f'(r)=O(r^{k-1})$ as $r \rightarrow %%@
0$. In particular,
\begin{equation}\label{zero}
\hbox{ $\varphi_{0}=1$ and  $f_{0}=\frac{1}{n \omega_{n}r^{n-1}}\int_{\partial B_{r}}u ds=
\frac{1}{n \omega_{n}}\int_{|x|=1}u(rx)ds.$}
\end{equation}
We also have  for any $k\geq 0$, and any continuous real valued functions $v$ and $w$ on %%@
$(0,\infty)$,
\begin{equation}
\int_{R^n}V(|x|)|\Delta u_{k}|^{2}dx=\int_{R^n}V(|x|)\big( \Delta %%@
f_{k}(|x|)-c_{k}\frac{f_{k}(|x|)}{|x|^2}\big)^{2}dx,
\end{equation}
and
\begin{equation}
\int_{R^n}W(|x|)|\nabla u_{k}|^{2}dx=\int_{R^n}W(|x|)|\nabla %%@
f_{k}|^{2}dx+c_{k}\int_{R^n}W(|x|)|x|^{-2}f^{2}_{k}dx.
\end{equation}

\begin{theorem} \label{main.hr} Let $V$ and $W$ be positive radial $C^1$-functions  on a ball %%@
$B\backslash \{0\}$, where $B$ is a ball with radius $R$ in $\R^n$ ($n \geq 1$) and centered at %%@
zero. Assume $\int^{R}_{0}\frac{1}{r^{n-1}V(r)}dr=\infty$ and $\lim_{r \rightarrow %%@
0}r^{\alpha}V(r)=0$ for some $\alpha<(n-2)$. If 
\begin{equation}\label{main.con}
W(r)-\frac{2V(r)}{r^2}+\frac{2V_r(r)}{r}-V_{rr}(r)\geq 0 \ \ for \ \ 0\leq r \leq R,
\end{equation} 
then the following statements are equivalent. 
\begin{enumerate}
\item $(V, W)$ is a Bessel pair with $\beta (V, W; R)\geq 1$.  
\item The  following inequality holds for all $u \in  
C^{\infty}_{0}(B)$,
\begin{equation*}\label{gen-hardy}
\hbox{   $({\rm HR}_{V,W})$ \quad \quad  $\int_{B}V(x)|\Delta u |^{2}dx \geq  \int_{B} W(x)|\nabla %%@
u|^{2}dx+(n-1)\int_{B}(\frac{V(x)}{|x|^2}-\frac{V_r(|x|)}{|x|})|\nabla u|^2dx.$  \quad \quad \quad %%@
\quad \quad \quad \quad  }
\end{equation*}

\end{enumerate} 
Moreover, if  $\beta (V, W; R)\geq 1$, then the best constant is given by 
\begin{equation}
\hbox{$\beta (V, W; R)=\sup\big\{c; \, \,  ({\rm HR}_{V, cW})$ holds$\big\}$.} 
\end{equation}

\end{theorem} 
{\bf Proof:}  That 2) implies 1) follows from Theorem \ref{mainrad.hr} and does not require %%@
condition (\ref{main.con}). 
To prove that 1) implies 2) assume that the equation  $(B_{V,W})$ has a positive solution on %%@
$(0,R]$. We prove that the inequality $(HR_{V,W})$ holds for all $u \in C^{\infty}_{0}(B)$ by %%@
frequently using that

\begin{equation}\label{1-dim}
\hbox{$\int^{R}_{0}V(r)|x'(r)|^2r^{n-1}dr \geq \int^{R}_{0}W(r)x^{2}(r)r^{n-1}dr$ for all $x\in %%@
C^1(0, R]$.}
\end{equation}

Indeed,  for all $n\geq 1$ and $k\geq 0$ we have 
\begin{eqnarray*}
\frac{1}{nw_n}\int_{R^n}V(x)|\Delta u_{k}|^{2}dx&=&\frac{1}{nw_n}\int_{R^n}V(x)\big( \Delta %%@
f_{k}(|x|)-c_{k}\frac{f_{k}(|x|)}{|x|^2}\big)^{2}dx\\
&=&
\int^{R}_{0}V(r)\big(f_{k}''(r)+\frac{n-1}{r}f_{k}'(r)-c_{k}\frac{f_{k}(r)}{r^2}\big)^{2}r^{n-1}dr
\\
&=&\int^{R}_{0}V(r)(f_{k}''(r))^{2}r^{n-1}dr+(n-1)^{2} \int^{R}_{0}V(r)(f_{k}'(r))^{2}r^{n-3}dr\\
&&+c^{2}_{k} \int^{R}_{0}V(r)f_{k}^{2}(r)r^{n-5}
+ 2(n-1) \int^{R}_{0}V(r)f_{k}''(r)f_{k}'(r)r^{n-2}\\
&&-2c_{k} 
\int^{R}_{0}V(r)f_{k}''(r)f_{k}(r)r^{n-3}dr
- 2c_{k}(n-1) \int^{R}_{0}V(r)f_{k}'(r)f_{k}(r)r^{n-4}dr.
\end{eqnarray*}
Integrate by parts and use (\ref{zero}) for $k=0$ to get
\begin{eqnarray}
\frac{1}{n\omega_{n}}\int_{R^n}V(x)|\Delta u_{k}|^{2}dx&=&  %%@
\int^{R}_{0}V(r)(f_{k}''(r))^{2}r^{n-1}dr+(n-1+2c_{k}) \int^{R}_{0}V(r)(f_{k}'(r))^{2}r^{n-3}dr %%@
\label{piece.0}\\
&+& 
(2c_{k}(n-4)+c^{2}_{k})\int^{R}_{0}V(r)r^{n-5}f_{k}^{2}(r)dr-(n-1)\int^{R}_{0}V_r(r)r^{n
-2}(f_{k}')^{2}(r)dr\nonumber\\
&-&c_{k}(n-5)\int^{R}_{0}V_r(r)f_{k}^2(r)r^{n-4}dr-c_{k}\int^{R}_{0}V_{rr}(r)f_{k}^2(r)r^{n-3}dr.
\nonumber
\end{eqnarray}
Now define $g_{k}(r)=\frac{f_{k}(r)}{r}$ and note that $g_{k}(r)=O(r^{k-1})$ for all $k\geq 1$. We %%@
have
\begin{eqnarray*}
\int^{R}_{0}V(r)(f_{k}'(r))^{2}r^{n-3}&=&\int^{R}_{0}V(r)(g_{k}'(r))^{2}r^{n-1}dr+\int^{R}_{0}2V(r
)g
_{k}(r)g_{k}'(r)r^{n-2}dr+\int^{R}_{0}V(r)g_{k}^{2}(r)r^{n-3}dr\\
&=&\int^{R}_{0}V(r)(g_{k}'(r))^{2}r^{n-1}dr-(n-3)\int^{R}_{0}V(r)g_{k}^{2}(r)r^{n-3}dr
-\int_{0}^{R}V_r(r)g^2_{k}(r)r^{n-2}dr\\
\end{eqnarray*}
Thus,
\begin{equation}\label{g1}
\int^{R}_{0}V(r)(f_{k}'(r))^{2}r^{n-3}\geq %%@
\int^{R}_{0}W(r)f_{k}^{2}(r)r^{n-3}dr-(n-3)\int^{R}_{0}V(r)f_{k}^{2}(r)r^{n-5}dr-\int_{0}^{R}V_r(r%%@
%%@
)f^2_{k}(r)r^{n-4}dr.
\end{equation}
Substituting $2c_k\int^{R}_{0}V(r)(f_{k}'(r))^{2}r^{n-3}$ in (\ref{piece.0}) by its lower estimate %%@
in the last inequality (\ref{g1}), we get
\begin{eqnarray*}
\frac{1}{n\omega_{n}}\int_{R^n}V(x)|\Delta u_{k}|^{2}dx&\geq&  
\int^{R}_{0}W(r)(f_{k}'(r))^{2}r^{n-1}dr+\int^{R}_{0}W(r)(f_{k}(r))^{2}r^{n-3}dr \\  &+&(n-1) %%@
\int^{R}_{0}V(r)(f_{k}'(r))^{2}r^{n-3}dr+c_{k}(n-1) \int^{R}_{0}V(r)(f_{k}(r))^{2}r^{n-5}dr\\
&-&(n-1)\int^{R}_{0}V_r(r)r^{n-2}(f_{k}')^{2}(r)dr-c_{k}(n-1)\int^{R}_{0}V_r(r
)r^{n
-4}(f_{k})^{2}(r)dr\\
&+&c_{k}(c_{k}-(n-1))\int^{R}_{0}V(r)r^{n-5}f_{k}^{2}(r)dr\\
&+&c_{k}\int^{R}_{0}(W(r)-\frac{2V(r)}{r^2}+\frac{2V_r(r)}{r}-V_{rr}(r))f^2_{k}(r)r^{n-3}dr
.\\
\end{eqnarray*}
The proof is now complete since the last term is non-negative by condition (\ref{main.con}). Note %%@
also that because of this condition, the formula for the best constant requires that $\beta (V, W;  %%@
R) \geq 1$, since if $W$ satisfies (\ref{main.con}) then $cW$ satisfies it for any $c\geq 1$.
\hfill $\square$
\begin{remark}\rm
In order to apply the above theorem to the Bessel pair 
\[
\hbox{$V(x)=|x|^{-2m}$ \quad and \quad $W_{m,c}(x)= 
V(x)\big[(\frac{n-2m-2}{2})^2|x|^{-2}+cW(x)\big]$}
\] 
where $W$ is a Bessel potential, we see that even in the simplest case $V\equiv 1$ and $W_{m,c}(x)= %%@
(\frac{n-2}{2})^2|x|^{-2}+W(x)$,  condition (\ref{main.con}) reduces to %%@
$(\frac{n-2}{2})^2|x|^{-2}+W(x) \geq 2|x|^{-2}$, which is then guaranteed only if $n\geq 5$. 

 More generally,  if $V(x)=|x|^{-2m}$, then in order to satisfy (\ref{main.con}) we need to have
\begin{equation}\label{restrict}
\frac{-(n+4)-2\sqrt{n^2-n+1}}{6} \leq m\leq \frac{-(n+4)+2\sqrt{n^2-n+1}}{6},  
\end{equation}
and in this case, we have for $m<\frac{n-2}{2}$ and any  Bessel potential $W$ on $B_{R}$, that  %%@
for all $u \in C^{\infty}_{0}(B_{R})$
\begin{equation}\label{gm-hr.0}
\int_{B_{R}}\frac{|\Delta u|^2}{|x|^{2m}}\geq  (\frac{n+2m}{2})^2\int_{B_{R}}\frac{|\nabla %%@
u|^2}{|x|^{2m+2}}dx+\beta (W; R)\int_{B_{R}}W(x)\frac{|\nabla u|^2}{|x|^{2m}}dx.
\end{equation}
Moreover, $(\frac{n+2m}{2})^2$ and $\beta (W; R)$ are the best constant. 

Therefore, inequality (\ref{gm-hr.0}) in the case where $m=0$ and $n\geq 5$,  already includes Theorem 1.5 in \cite{TZ} as a special case. It also extends Theorem 1.8 in \cite{TZ} where it is  established under 
the condition 
\begin{equation}
0\leq m\leq \frac{-(n+4)+2\sqrt{n^2-n+1}}{6}
\end{equation} which is more restrictive than (\ref{restrict}). 
We shall see however that %%@
this inequality remains true without condition (\ref{restrict}), but with a  constant 
that is sometimes different from $(\frac{n+2m}{2})^2$ in the cases where (\ref{restrict}) is not %%@
valid. For example, if $m=0$, then the best constant is $3$ in dimension $4$ and $\frac{25}{36}$ %%@
in dimension $3$. 
\end{remark}
%\begin{remark} \rm In view of the above remarks, Theorem \ref{main.hr} in the case where $m=0$ and %%@
%$n\geq 5$ already includes Theorem 1.5 in \cite{TZ} as a special case. Moreover,  when %%@
%$V(x)=|x|^{-2m}$,  the above shows that Theorem 1.8 in \cite{TZ} is still valid if one replaces %%@
%the condition 
%\[0\leq m\leq \frac{-(n+4)+2\sqrt{n^2-n+1}}{6},\]
%by
%\[\frac{-(n+4)-2\sqrt{n^2-n+1}}{6} \leq m\leq \frac{-(n+4)+2\sqrt{n^2-n+1}}{6}.
%\]
%\end{remark}
We shall now give a few immediate applications of the above in the case where $m=0$ and $n\geq 5$.
Actually the results are true in lower dimensions, and will be stated as such,  but the proofs for %%@
$n<5$ will require additional work and will be postponed to the next section.  

%and $n\geq 3$. 
 
\begin{theorem} \label{m=0}Assume $W$ is a Bessel potential on $B_{R} \subset R^n$ with $n\geq 3$,  %%@
then for all $u \in C_{0}^{\infty}(B_{R})$ we have 
\begin{equation}
\int_{B_{R}}|\Delta u|^2 dx\geq C(n)\int_{B_{R}}\frac{|\nabla u|^2}{|x|^2}dx+\beta(W;  %%@
R)\int_{B_{R}}W(x)|\nabla u|^2dx,
\end{equation} 
where $C(3)=\frac{25}{36}$, $C(4)=3$ and $C(n)=\frac{n^2}{4}$ for all $n\geq 5$. Moreover, 
 $C(n)$ and $\beta (W;  R)$ are the best constants. 
 
In particular,  the following holds for any smooth bounded domain $\Omega$ in $R^n$ with %%@
$R=\sup_{x \in \Omega}|x|$, and any $u \in H^{2}(\Omega)\cap H^{1}_{0}(\Omega)$.
\begin{itemize}
\item  For any $\alpha<2$, 
\begin{equation}
\int_{\Omega}|\Delta u|^2 dx\geq C(n)\int_{\Omega}\frac{|\nabla u|^2}{|x|^2}dx+\beta(|x|^\alpha;  %%@
R)\int_{\Omega}\frac{|\nabla u|^2}{|x|^\alpha}dx,
\end{equation} 
and for $\alpha=0$,
 \begin{equation}\label{gm-hardy-Rellich}
\hbox{$\int_{\Omega}|\Delta u |^{2}dx \geq C(n) \int_{\Omega}\frac{|\nabla %%@
u|^2}{|x|^{2}}dx+\frac{z^{2}_{0}}{R^2}\int_{\Omega}|\nabla u|^{2}dx$,}
\end{equation}
the constants being optimal when $\Omega$ is a ball. 
\item For any $k\geq 1$, and $\rho=R( e^{e^{e^{.^{.^{e(k-times)}}}}} )$, we have 
 \begin{equation}
 \int_{\Omega}|\Delta u(x) |^{2}dx \geq C(n)\int_{\Omega}\frac{|\nabla u|^2}{|x|^2} %%@
dx+\frac{1}{4}\sum^{k}_{j=1}\int_{\Omega}\frac{|\nabla u |^2}{|x|^2}\big( %%@
\prod^{j}_{i=1}log^{(i)}\frac{\rho}{|x|}\big)^{-2}dx, 
 \end{equation} 
 \item For $D\geq R$, and $X_{i}$ is defined as (\ref{X-def}) we have
 \begin{equation}
  \int_{\Omega}|\Delta u(x) |^{2}dx \geq C(n)\int_{\Omega}\frac{|\nabla u|^2}{|x|^2} %%@
dx+\frac{1}{4}\sum^{\infty}_{i=1}\int_{\Omega}\frac{|\nabla %%@
u|}{|x|^{2}}X^{2}_{1}(\frac{|x|}{D})X^{2}_{2
}(\frac{|x|}{D})...X^{2}_{i}(\frac{|x|}{D})dx,  
\end{equation} 
 Moreover, all constants appearing in the above two
inequality are optimal.
\end{itemize}
\end{theorem} 

\begin{theorem} Let $W(x)=W(|x|)$ be radial Bessel potential on a ball $B$ of radius $R$ in $R^n$ %%@
with $n\geq 4$, and such that $\frac{W_r(r)}{W(r)}=\frac{\lambda}{r}+f(r)$, where $f(r)\geq 0$ and %%@
$\lim_{r \rightarrow 0}rf(r)=0$.  If $\lambda <n-2$,  then the following Hardy-Rellich inequality %%@
holds:
\begin{equation}
\int_{B}|\Delta u|^{2}dx \geq %%@
\frac{n^{2}(n-4)^2}{16}\int_{B}\frac{u^2}{|x|^4}dx+(\frac{n^2}{4}+\frac{(n-\lambda-2)^2}{4})
\beta (W;  R)\int_{B}\frac{W(x)}{|x|^2}u^2 dx,
\end{equation}
\end{theorem}
{\bf Proof:} Use first Theorem \ref{m=0} with the Bessel potential $W$, then Theorem %%@
\ref{main-CKN} with the Bessel pair\\ 
$(|x|^{-2}, |x|^{-2}(\frac{(n-4)^2}{4}|x|^{-2}+W)$, then Theorem \ref{super.hardy} with the Bessel %%@
pair $(W, \frac{(n-\lambda-2)^2}{4})|x|^{-2}W)$ to obtain 
\begin{eqnarray*}
\int_{B}|\Delta u|^2 dx&\geq&C(n)\int_{B}\frac{|\nabla u|^2}{|x|^2}dx+ \beta (W, %%@
R)\int_{B}W(x)|\nabla u|^2 dx\\
&\geq& C(n)\frac{(n-4)^2}{4}\int_{B}\frac{u^2}{|x|^4}dx+C(n)\beta (W, %%@
R)\int_{B}\frac{W(x)}{|x|^2}u^2+\beta (W, R)\int W(x)|\nabla u|^2 dx\\
&\geq&
C(n)\frac{ (n-4)^2}{4}\int_{B}\frac{u^2}{|x|^4}dx+(C(n)+\frac{(n-\lambda-2)^2}{4})
\beta (W, R)\int_{B}\frac{W(x)}{|x|^2}u^2 dx.
\end{eqnarray*}
Recall that $C(n)=\frac{n^2}{4}$ for $n\geq 5$, giving the claimed result in these dimensions. %%@
This is however not the case when $n=4$, and therefore another proof will be given in the next %%@
section to cover these cases. 

The following is immediate from Theorem \ref{m=0} and from the fact that $\lambda=2$ for the %%@
Bessel potential under consideration. 
 
\begin{corollary} Let $\Omega$ be a smooth bounded domain in $\R^n$, $n \geq 4$ and $R=\sup_{x \in %%@
\Omega}|x|$. Then the following holds for all $u \in H^{2}(\Omega) \cap H^{1}_{0}(\Omega)$
\begin{enumerate}
\item  If $\rho=R( e^{e^{e^{.^{.^{e(k-times)}}}}} )$ and
$log^{(i)}(.)$ is defined as (\ref{log-def}), then 
\begin{equation}
 \int_{\Omega}|\Delta u(x) |^{2}dx \geq \frac{n^2(n-4)^2}{16}\int_{\Omega}\frac{u^2}{|x|^4} %%@
dx+(1+\frac{n(n-4)}{8})\sum^{k}_{j=1}\int_{\Omega}\frac{u^2}{|x|^4}\big( %%@
\prod^{j}_{i=1}log^{(i)}\frac{\rho}{|x|}\big)^{-2}dx.
 \end{equation} 
\item If $D\geq R$ and $X_{i}$ is defined as (\ref{X-def}), then 
\begin{equation}
\int_{\Omega}|\Delta u(x) |^{2}dx \geq \frac{n^2(n-4)^2}{16}\int_{\Omega}\frac{u^2}{|x|^4} %%@
dx+(1+\frac{n(n-4)}{8})\sum^{\infty}_{i=1}\int_{\Omega}\frac{u^2}{|x|^{4}}X^{2}_{1}(\frac{|x|}{D})%%@
X^{2}_{2
}(\frac{|x|}{D})...X^{2}_{i}(\frac{|x|}{D})dx.
\end{equation} 
\end{enumerate}
  \end{corollary} 

%We will prove the above two theorems later in this section.

\begin{theorem} \label{n-dim} Let $W_1(x)$ and $W_2(x)$ be two radial Bessel potentials on a ball %%@
$B$ of radius $R$ in $R^n$ with $n\geq 4$. If $a<1$,  then  there exists $c(a, R)>0$ such that for %%@
all  $u \in H^{2}(B)\cap H_{0}^{1}(B)$
% \begin{eqnarray*}
%   \int_{\Omega}|\Delta u |^{2}dx&\geq&
%\frac{ C(n)(n-4)^2}{4}\int_{B}\frac{u^2}{|x|^4}dx+C(n)\beta (W_1; %%@
%R)\int_{B}W_1(x)\frac{u^2}{|x|^2}dx\\
%&&+ c(\frac{n-2a-2}{2})^2 \int_{\Omega}\frac{u^2}{|x|^{2a+2}}dx
%+c\beta (W_2; R)\int_{\Omega}W_2(x)\frac{u^2}{|x|^{2a}}dx.
% \end{eqnarray*}
%In particular, if $n\geq 5$
 \begin{eqnarray*}
   \int_{B}|\Delta u |^{2}dx &\geq& \frac{n^2(n-4)^2}{16} %%@
\int_{B}\frac{u^{2}}{|x|^{4}}dx+ \frac{n^2}{4}\beta (W_1;  R)\int_{B} %%@
W_1(x)\frac{u^{2}}{|x|^2}dx\\
&&+ c(\frac{n-2a-2}{2})^2 \int_{B}\frac{u^2}{|x|^{2a+2}}dx+c\beta (W_2;  %%@
R)\int_{B}W_2(x)\frac{u^2}{|x|^{2a}}dx, 
\end{eqnarray*}
% and if $n=4$, then
 % \begin{eqnarray*}
 %  \int_{B}|\Delta u |^{2}dx&\geq&
%3\beta (W_1;  R)\int_{B}W_1(x)\frac{u^2}{|x|^2}dx+ c(1-a)^2 \int_{B}\frac{u^2}{|x|^{2a+2}}dx
%+c\beta (W_2;  R)\int_{B}W_2(x)\frac{u^2}{|x|^{2a}}dx.
% \end{eqnarray*}

\end{theorem} 
{\bf Proof:} Here again we shall give a proof when $n\geq 5$. The case $n=4$ will be handled in %%@
the next section. We again  first use Theorem \ref{m=0} (for $n\geq 5$) with the Bessel potential %%@
$|x|^{-2a}$ where $a<1$, then Theorem \ref{main-CKN} with the Bessel pair  $(|x|^{-2}, %%@
|x|^{-2}(\frac{(n-4)^2}{4}|x|^{-2}+W) )$, then again Theorem \ref{main-CKN} with the Bessel pair\\
 $(|x|^{-2a}, |x|^{-2a}((\frac{n-2a-2}{2})^2|x|^{-2}+W)$
 to obtain 
\begin{eqnarray*}
\int_{B}|\Delta u|^2 dx&\geq& \frac{n^2}{4} \int_{B}\frac{|\nabla u|^2}{|x|^2}dx+ \beta %%@
(|x|^{-2a};  R)\int_{B}\frac{|\nabla u|^2}{|x|^{-2a}} dx\\
&\geq&  \frac{n^2(n-4)^2}{16}\int_{B}\frac{u^2}{|x|^4}dx+\frac{n^2}{4} \beta (W_1;  %%@
R)\int_{B}W_1(x)\frac{u^2}{|x|^2}dx+ \beta (|x|^{-2a};  R)\int_{B}\frac{|\nabla u|^2}{|x|^{-2a}} %%@
dx\\&\geq&
\frac{n^2(n-4)^2}{16}\int_{B}\frac{u^2}{|x|^4}dx+\frac{n^2}{4} \beta (W_1;  %%@
R)\int_{B}W_1(x)\frac{u^2}{|x|^2}dx\\
&&+ \beta (|x|^{-2a};  R)(\frac{n-2a-2}{2})^2 \int_{B}\frac{u^2}{|x|^{2a+2}}dx
+\beta (|x|^{-2a}; R)\beta(W_2;  R)\int_{B}W_2(x)\frac{u^2}{|x|^{2a}}dx.
 \end{eqnarray*}
The following theorem will be established in full generality (i.e with $V(r)=r^{-m}$) in the next %%@
section.

\begin{theorem} Let $W(x)=W(|x|)$ be a radial Bessel potential on a smooth bounded domain $\Omega$ %%@
in $\R^n$, $n \geq 4$. Then, 
\begin{equation*}\label{HR1}
\hbox{$   \quad \quad \quad  \int_{\Omega}|\Delta u(x) |^{2}dx -\frac{n^2(n-4)^2}{16} %%@
\int_{\Omega}\frac{u^{2}}{|x|^{4}}dx-
\frac{n^2}{4}\int_{\Omega} W(x)u^{2}dx\geq \frac{z^{2}_{0}}{R^2}||u||^{2}_{W_{0}^{1,2}(\Omega)} %%@
\quad u \in H^{2}_{0}(\Omega). \quad \quad $}
\end{equation*} 
\end{theorem}

\subsection{The case of  power potentials $|x|^m$}

The general Theorem \ref{main.hr} allowed us to deduce inequality (\ref{gm-hr}) below for a %%@
restricted interval of powers $m$. We shall now prove that the same holds for all  %%@
$m<\frac{n-2}{2}$.  The following theorem improves considerably Theorem 1.7, Theorem 1.8, and %%@
Theorem 6.4 in \cite{TZ}.
 \begin{theorem}\label{gm.hr}
Suppose $n\geq 1$ and $m<\frac{n-2}{2}$, and let $W$ be a Bessel potential on a ball $B_{R}\subset %%@
R^n$ of radius $R$. Then for all $u \in C^{\infty}_{0}(B_{R})$
\begin{equation}\label{gm-hr}
\int_{B_{R}}\frac{|\Delta u|^2}{|x|^{2m}}\geq a_{n,m}\int_{B_{R}}\frac{|\nabla %%@
u|^2}{|x|^{2m+2}}dx+\beta(W;  R)\int_{B_{R}}W(x)\frac{|\nabla u|^2}{|x|^{2m}}dx,
\end{equation}
where 
\[
a_{n,m}=\inf \left\{\frac{\int_{B_{R}}\frac{|\Delta u|^2}{|x|^{2m}}dx}{\int_{B_{R}}\frac{|\nabla %%@
u|^2}{|x|^{2m+2}}dx};\, u\in C^{\infty}_{0}(B_{R})\setminus \{0\}\right\}.
\]
Moreover, $\beta(W;  R)$ and $a_{m,n}$ are the best constants to be computed in the appendix. 
\end{theorem}
{\bf Proof:}  Assuming  the inequality
\[\int_{B_{R}}\frac{|\Delta u|^2}{|x|^{2m}}\geq a_{n,m}\int_{B_{R}}\frac{|\nabla %%@
u|^2}{|x|^{2m+2}}dx,\]
holds for all $u \in C^{\infty}_{0}(B_{R})$, we shall prove that it can be improved by any Bessel %%@
potential $W$. We will use the following inequality frequently in the proof which follows directly %%@
from Theorem \ref{main-CKN} with n=1.
\begin{equation}\label{freq-in}
\int_{0}^{R}r^{\alpha}(f'(r))^{2}dr\geq %%@
(\frac{\alpha-1}{2})^2\int_{0}^{R}r^{\alpha-2}f^2(r)dr+\beta(W;R)\int_{0}^{R}r^{\alpha}W(r)
f^2(r)dr, \ \ \alpha \geq 1,
\end{equation}
for all $f \in C^{\infty}(0,R)$, where both $(\frac{\alpha-1}{2})^2$ and $\beta(W;R)$ are best %%@
constants. 

Decompose $u \in C^{\infty}_{0}(B_{R})$ into its spherical harmonics $ %%@
\Sigma^{\infty}_{k=0}u_{k}$, where $u_{k}=f_{k}(|x|)\varphi_{k}(x)$. We evaluate %%@
$I_k=\frac{1}{nw_n}\int_{R^n} \frac{|\Delta u_{k}|^2}{|x|^{2m}}dx$ in the following way
\begin{eqnarray*}
 I_k&=&\int_{0}^{R}r^{n-2m-1}(f''_{k}(r))^2dr+[(n-1)(2m+1)+2c_{k}]\int^{R}_{0}r^{n-2m-3}(f_{k}')^2
dr\\
&&+c_{k}[c_{k}+(n-2m-4)(2m+2)]\int_{0}^{R}r^{n-2m-5}(f_{k}(r))^2dr\\
&\geq& \beta %%@
(W)\int^{R}_{0}r^{n-2m-1}W(x)(f_{k}')^2dr+[(\frac{n+2m}{2})^2+2c_{k}]\int^{R}_{0}r^{n-2m-3}(f_{k}'%%@
)^2dr\\
&&+c_{k}[c_{k}+(n-2m-4)(2m+2)]\int_{0}^{R}r^{n-2m-5}(f_{k}(r))^2dr\\
&\geq& \beta (W)\int^{R}_{0}r^{n-2m-1}W(x)(f_{k}')^2dr+a_{n,m}\int^{R}_{0}r^{n-2m-3}(f_{k}')^2dr\\
&&+\beta (W)[(\frac{n+2m}{2})^2+2c_{k}-a_{n,m}]\int^{R}_{0}r^{n-2m-3}W(x)(f_{k})^2dr\\
&&+\big((\frac{n-2m-4}{2})^2[(\frac{n+2m}{2})^2+2c_{k}-a_{n,m}]+c_{k}[c_{k}+(n-2m-4)(2m+2)]\big)
\int_{0}^{R}r^{n-2m-5}(f_{k}(r))^2dr.\\
\end{eqnarray*}
Now by (\ref{a-nm-k}) we have 
\[\big((\frac{n-2m-4}{2})^2[(\frac{n+2m}{2})^2+2c_{k}-a_{n,m}]+c_{k}[c_{k}+(n-2m-4)(2m+2)]\geq %%@
c_{k}a_{n,m},\]
for all $k\geq 0$. Hence, we have
\begin{eqnarray*}
I_{k}&\geq&a_{n,m}\int^{R}_{0}r^{n-2m-3}(f_{k}')^2dr+a_{n,m}c_{k}\int_{0}^{R}r^{n-2m-5}(f_{k}(r))^%%@
2dr\\
&&+\beta (W)\int^{R}_{0}r^{n-2m-1}W(x)(f_{k}')^2dr+\beta %%@
(W)[(\frac{n+2m}{2})^2+2c_{k}-a_{n,m}]\int^{R}_{0}r
^{n-2m-3}W(x)(f_{k})^2dr\\
&\geq& a_{n,m}\int^{R}_{0}r^{n-2m-3}(f_{k}')^2dr+a_{n,m}c_{k} %%@
\int_{0}^{R}r^{n-2m-5}(f_{k}(r))^2dr\\
&&+ \beta (W)\int^{R}_{0}r^{n-2m-1}W(x)(f_{k}')^2dr+
\beta (W)c_{k}\int^{R}_{0}r^{n-2m-3}W(x)(f_{k})^2dr\\
&=&a_{n,m}\int_{B_{R}}\frac{|\nabla u|^2}{|x|^{2m+2}}dx+\beta (W)\int_{B_{R}}W(x)\frac{|\nabla %%@
u|^2}{|x|^{2m}}dx.
\end{eqnarray*}
Moreover, it is easy to see from Theorem \ref{main} and the above calculation that $\beta (W; R)$ %%@
is the best constant.   

\begin{theorem}\label{super.hardy-rellich}
Let $\Omega$ be a smooth domain in $R^{n}$ with $n\geq 1$ and let $V \in C^{2}(0,R=:\sup_{x \in %%@
\Omega}|x|)$ be a non-negative function that satisfies the following conditions:
\begin{equation}
\hbox{$V_r(r)\leq 0$\quad and \quad 
$\int^{R}_{0}\frac{1}{r^{n-3}V(r)}dr=-\int^{R}_{0}\frac{1}{r^{n-4}V_r(r)}dr=+\infty$.}
\end{equation}
 There exists $\lambda_{1}, \lambda_{2} \in R$ such that 
\begin{equation}
\hbox{$\frac{rV_r(r)}{V(r)}+\lambda_{1} \geq 0$ on $(0, R)$ and $\lim\limits_{r\to %%@
0}\frac{rV_r(r)}{V(r)}+\lambda_{1} =0$,}
\end{equation}
\begin{equation}
\hbox{$\frac{rV_{rr}(r)}{V_r(r)}+\lambda_{2} \geq 0$ on $(0, R)$ and $\lim\limits_{r\to %%@
0}\frac{rV_{rr}(r)}{V_r(r)}+\lambda_{2} =0$,}
\end{equation}
and
\begin{equation}\label{super.hr.con}
\hbox{$\left(\frac{1}{2}(n-\lambda_{1}-2)^2+3(n-3)\right)V(r)-(n-5)rV_r(r)-r^2V_{rr}(r)\geq 0$ for %%@
all $r \in (0,R)$. }
\end{equation}
Then the following inequality holds:
\begin{eqnarray}\label{super.hr}
\int_{\Omega}V(|x|)|\Delta u|^2 dx&\geq& %%@
(\frac{(n-\lambda_{1}-2)^{2}}{4}+(n-1))\frac{(n-\lambda_{1}-4)^{2}}{4}\int_{\Omega}\frac{V(|x|)}{|%%@
x|^4}u^2 dx \nonumber\\
&&-\frac{(n-1)(n-\lambda_{2}-2)^{2}}{4}\int_{\Omega}\frac{V_r(|x|)}{|x|^3}u^2 dx.
\end{eqnarray}
\end{theorem}
{\bf Proof:} We have by Theorem \ref{super.hardy} and condition (\ref{super.hr.con}), 
{\small
\begin{eqnarray*}
\frac{1}{n\omega_{n}}\int_{R^n}V(x)|\Delta u_{k}|^{2}dx&=&  %%@
\int^{R}_{0}V(r)(f_{k}''(r))^{2}r^{n-1}dr+(n-1+2c_{k}) \int^{R}_{0}V(r)(f_{k}'(r))^{2}r^{n-3}dr %%@
\label{piece}\\
&+& 
(2c_{k}(n-4)+c^{2}_{k})\int^{R}_{0}V(r)r^{n-5}f_{k}^{2}(r)dr-(n-1)\int^{R}_{0}V_r(r)r^{n
-2}(f_{k}')^{2}(r)dr\\
&-&c_{k}(n-5)\int^{R}_{0}V_r(r)f_{k}^2(r)r^{n-4}dr-c_{k}\int^{R}_{0}V_{rr}(r)f_{k}^2(r)r^{n-3}dr\\
&\geq&\int^{R}_{0}V(r)(f_{k}''(r))^{2}r^{n-1}dr+(n-1) \int^{R}_{0}V(r)(f_{k}'(r))^{2}r^{n-3}dr %%@
\label{piece}\\
&-& 
(n-1)\int^{R}_{0}V_r(r)r^{n-2}(f_{k}')^{2}(r)dr\\
&+&c_{k}\int_{0}^{R}\left(\left(\frac{1}{2}(n-\lambda_{1}-2)^2+3(n-3)\right)V(r)-(n-5)rV_r(r)-r^2V%%@
%%@
_{rr}(r) \right)f_{k}^2(r)r^{n-5}dr
%&\times& f_{k}^2(r)r^{n-5}dr,
\end{eqnarray*}
}
The rest of the proof follows %%@
from the above inequality combined with Theorem \ref{super.hardy}. \hfill $\Box$
\begin{remark}\rm
Let $V(r)=r^{-2m}$ with $m\leq \frac{n-4}{2}$. Then in order to satisfy condition %%@
(\ref{super.hr.con}) we must have $-1-\frac{\sqrt{1+(n-1)^2}}{2}\leq m\leq \frac{n-4}{2}$. Under %%@
this assumption the inequality (\ref{super.hr}) gives the following weighted second order Rellich %%@
inequality: 
\[\int_{B}\frac{|\Delta u|^2}{|x|^{2m}}dx\geq %%@
(\frac{(n+2m)(n-4-2m)}{4})^2\int_{B}\frac{u^2}{|x|^{2m+4}}dx.\]
In the following theorem we will show that the constant appearing in the above inequality is %%@
optimal. Moreover, we will see that if $m<-1-\frac{\sqrt{1+(n-1)^2}}{2}$, then the best constant %%@
is strictly less than $(\frac{(n+2m)(n-4-2m)}{4})^2$. This shows that inequality (\ref{super.hr}) %%@
is actually sharp.
\end{remark}

\begin{theorem}
Let $m\leq \frac{n-4}{2}$ and define
\begin{equation}
\beta_{n,m}=\inf_{u \in C^{\infty}_{0}(B)\backslash \{0\}}\frac{\int_{B}\frac{|\Delta %%@
u|^2}{|x|^{2m}}dx}{\int_{B}\frac{u^2}{|x|^{2m+4}}dx}.
\end{equation}
Then
\[
\beta_{n,m}=(\frac{(n+2m)(n-4-2m)}{4})^2+\min_{k=0,1,2,...}\{
k(n+k-2)[k(n+k-2)+\frac{(n+2m)(n-2m-4)}{2}]\}.\]
Consequently the values of $\beta_{n,m}$ are as follows.
\begin{enumerate}
\item If $-1-\frac{\sqrt{1+(n-1)^2}}{2}\leq m\leq \frac{n-4}{2}$, then
\[\beta_{n,m}=(\frac{(n+2m)(n-4-2m)}{4})^2.\]
\item If $\frac{n}{2}-3 \leq m\leq -1-\frac{\sqrt{1+(n-1)^2}}{2}$, then
\[\beta_{n,m}=(\frac{(n+2m)(n-4-2m)}{4})^2+(n-1)[(n-1)+\frac{(n+2m)(n-2m-4)}{2}].\]
\item If $k:=\frac{n-2m-4}{2} \in N$, then
\[\beta_{n,m}=(\frac{(n+2m)(n-4-2m)}{4})^2+k(n+k-2)[k(n+k-2)+\frac{(n+2m)(n-2m-4)}{2}].\]
\item If $k<\frac{n-2m-4}{2}<k+1$ for some $k \in N$, then
\begin{eqnarray*}
\beta_{n,m}=\frac{(n+2m)^2(n-2m-4)^2}{16}+a(m, n, k)
\end{eqnarray*}
\end{enumerate}
where 
\footnotesize
\[
a(m,n,k)=\min\left\{k(n+k-2)[k(n+k-2)+\frac{(n+2m)(n-2m-4)}{2}], %%@
(k+1)(n+k-1)[(k+1)(n+k-1)+\frac{(n+2m)(n-2m-4)}{2}]\right\}.
\]

\end{theorem}
{\bf Proof:} Decompose $u \in C^{\infty}_{0}(B_{R})$ into spherical harmonics 
$\Sigma^{\infty}_{k=0}u_{k}$, where $u_{k}=f_{k}(|x|)\varphi_{k}(x)$. we have
\begin{eqnarray*}
\frac{1}{n\omega_{n}}\int_{R^n} \frac{|\Delta %%@
u_{k}|^2}{|x|^{2m}}dx&=&\int_{0}^{R}r^{n-2m-1}(f''_{k}(r))^2dr+[(n-1)(2m+1)+2c_{k}]\int^{R}_{0}r^{%%@
n-2m-3}(f_{k}')^2dr\\
&+&c_{k}[c_{k}+(n-2m-4)(2m+2)]\int_{0}^{R}r^{n-2m-5}(f_{k}(r))^2dr\\
&\geq&\big((\frac{(n+2m)(n-4-2m)}{4})^2\\
&+&c_{k}[c_{k}+\frac{(n+2m)(n-2m-4)}{2}]\big)\int_{0}^{R}r^{n-2m-5}(f_{k}(r))^2dr,
\end{eqnarray*}
by Hardy inequality. Hence,
\[\beta_{n,m}\geq B(n,m,k):= (\frac{(n+2m)(n-4-2m)}{4})^2
+ \min_{k=0,1,2,...}\{k(n+k-2)[k(n+k-2)+\frac{(n+2m)(n-2m-4)}{2}]\}.\]
To prove that $\beta_{n,m}$ is the best constant, let $k$ be such that 
\begin{eqnarray}
\beta_{n,m}=\frac{(n+2m)(n-4-2m)}{4})^2
+k(n+k-2)[k(n+k-2)+\frac{(n+2m)(n-2m-4)}{2}].
\end{eqnarray}
Set 
\[u=|x|^{-\frac{n-4}{2}+m+\epsilon}\varphi_{k}(x)\varphi(|x|),\]
where $\varphi_{k}(x)$ is an eigenfunction corresponding to the eigenvalue $c_{k}$ and %%@
$\varphi(r)$ is a smooth cutoff function, such that $0 \leq \varphi \leq 1$, with $\varphi\equiv %%@
1$ in $[0,\frac{1}{2}]$. We have
\[\frac{\int_{B_R}\frac{|\Delta u|^2}{|x|^{2m}}dx}{\int_{B_{R}}\frac{ %%@
u^2}{|x|^{2m+4}}dx}=(-\frac{(n+2m)(n-4-2m)}{4}-c_{k}+\epsilon(2+2m+\epsilon))^2+O(1).\]
Let now $\epsilon \rightarrow 0$ to obtain the result. Thus the inequality
\[\int_{B_R}\frac{|\Delta u|^2}{|x|^{2m}}\geq \beta_{n,m}\int_{B_R}\frac{ %%@
u^2}{|x|^{2m+4}}dx,\]
holds for all $u \in C^{\infty}_{0}(B_{R})$. 

To calculate explicit values of $\beta_{n,m}$ we need to find the minimum point of the function
\[f(x)=x(x+\frac{(n+2m)(n-2m-4)}{2}), \ \ x\geq 0.\]
Observe that
\[f'(-\frac{(n+2m)(n-2m-4)}{4})=0.\]
To find minimizer $k \in N$ we should solve the equation
\[k^2+(n-2)k+\frac{(n+2m)(n-2m-4)}{4}=0.\]
The roots of the above equation are $x_{1}=\frac{n+2m}{2}$ and $x_{2}=\frac{n-2m-4}{2}$. 1) %%@
follows from Theorem \ref{super.hardy-rellich}. It is easy to see that if $m\leq %%@
-1-\frac{\sqrt{1+(n-1)^2}}{2}$, then $x_{1}<0$. Hence, for $m\leq -1-\frac{\sqrt{1+(n-1)^2}}{2}$ %%@
the minimum of the function $f$ is attained in $x_{2}$. Note that if $m\leq %%@
-1-\frac{\sqrt{1+(n-1)^2}}{2}$, then $B(n,m1)\leq B(n,m,0)$. Therefore claims 2), 3), and 4) %%@
follow.  \hfill $\Box$\\

The following theorem extends Theorem 1.6 of \cite{TZ} in many ways. First,  we do not assume that %%@
$n\geq 5$ or $m\geq 0$,  as was assumed there. Moreover, inequality (\ref{ex-gen-hr}) below  %%@
includes inequalities (1.17) and (1.22) of \cite{TZ} as special cases.

\begin{theorem}\label{gen.hr} Let $m\leq \frac{n-4}{2}$ and let $W(x)$ be a Bessel potential on a %%@
ball $B$ of radius $R$  in $R^n$ with radius $R$. Assume %%@
$\frac{W(r)}{W_r(r)}=-\frac{\lambda}{r}+f(r)$, where $f(r)\geq 0$ and $\lim_{r \rightarrow %%@
0}rf(r)=0$. Then the following inequality holds for all $u \in C^{\infty}_{0}(B)$
\begin{eqnarray}\label{ex-gen-hr}
\int_{B}\frac{|\Delta u|^{2}}{|x|^{2m}}dx &\geq& \beta_{n,m}\int_{B}\frac{u^2}{|x|^{2m+4}}dx %%@
\nonumber\\
&&\quad+\beta (W;  R)(\frac{(n+2m)^2}{4}+\frac{(n-2
m-\lambda-2)^2}{4})
\int_{B}\frac{W(x)}{|x|^{2m+2}}u^2 dx.
\end{eqnarray}
\end{theorem}

{\bf Proof:} Again we will frequently use inequality (\ref{freq-in})  in the proof. Decomposing $u %%@
\in C^{\infty}_{0}(B_{R})$ into spherical harmonics 
$\Sigma^{\infty}_{k=0}u_{k}$, where $u_{k}=f_{k}(|x|)\varphi_{k}(x)$, we can write
\begin{eqnarray*}
\frac{1}{n\omega_{n}}\int_{R^n} \frac{|\Delta %%@
u_{k}|^2}{|x|^{2m}}dx&=&\int_{0}^{R}r^{n-2m-1}(f''_{k}(r))^2dr+[(n-1)(2m+1)+2c_{k}]\int^{R}_{0}r^{%%@
n-2m-3}(f_{k}')^2dr\\
&&+c_{k}[c_{k}+(n-2m-4)(2m+2)]\int_{0}^{R}r^{n-2m-5}(f_{k}(r))^2dr\\
&\geq&(\frac{n+2m}{2})^2\int^{R}_{0}r^{n-2m-3}(f_{k}')^2dr+ \beta (W;  %%@
R)\int^{R}_{0}r^{n-2m-1}W(x)(f_{k}')^2dr\\
&&+c_{k}[c_{k}+2(\frac{n-\lambda-4}{2})^2+(n-2m-4)(2m+2)]\int_{0}^{R}r^{n-2m-5}(f_{k}(r))^2dr,
\end{eqnarray*}
where we have used the fact that $c_k\geq 0$ to get the above inequality. We have
\begin{eqnarray*}
\frac{1}{n\omega_{n}}\int_{R^n} \frac{|\Delta u_{k}|^2}{|x|^{2m}}dx
&\geq& \beta_{n,m} \int^{R}_{0}r^{n-2m-5}(f_{k})^2dr\\
&&+\beta (W;  R)\frac{(n+2m)^{2}}{4}
\int^{R}_{0}r^{n-2m-3}W(x)(f_{k})^2dr\\
&& +\beta (W;  R)\int^{R}_{0}r^{n-2m-1}W(x)(f_{k}')^2dr\\
&\geq&\beta_{n,m}\int^{R}_{0}r^{n-2m-5}(f_{k})^2dr\\
&&+\beta (W;  R)(\frac{(n+2m)^2}{4}+\frac{(n-2
m-\lambda-2)^2}{4})
\int^{R}_{0}r^{n-2m-3}W(x)(f_{k})^2dr \\
&\geq& \frac{\beta_{n,m}}{n\omega_{n}}\int_{B}\frac{u_{k}^2}{|x|^{2m+4}}dx\\
&&+\frac{\beta (W;  R)}{n\omega_{n}}(\frac{(n+2m)^2}{4}+\frac{(n-2
m-\lambda-2)^2}{4})
\int_{B}\frac{W(x)}{|x|^{2m+2}}u_{k}^2 dx,
\end{eqnarray*}
by Theorem \ref{super.hardy}. Hence, (\ref{ex-gen-hr}) holds and the proof is complete. \hfill %%@
$\Box$

\begin{theorem}\label{hrs-in} Assume $-1<m\leq \frac{n-4}{2}$ and let $W(x)$ be a Bessel potential %%@
on a ball $B$ of radius $R$ and centered at zero in $R^n$ ($n\geq 1$). Then there exists $C>0$ %%@
such that the following holds for all $u \in C^{\infty}_{0}(B)$:
\begin{eqnarray}\label{hrs}
\int_{B}\frac{|\Delta u|^{2}}{|x|^{2m}}dx &\geq&
\frac{(n+2m)^{2}(n-2m-4)^2}{16}\int_{B}\frac{u^2}{|x|^{2m+4}}dx  \\
&&+\beta (W; R)\frac{(n+2m)^2}{4}
\int_{B}\frac{W(x)}{|x|^{2m+2}}u^2 dx+ \beta(|x|^{2m}; R)||u||_{H^1_0}.
\end{eqnarray}
\end{theorem}
{\bf Proof:} Decomposing again $u \in C^{\infty}_{0}(B_{R})$ into its spherical harmonics %%@
$\Sigma^{\infty}_{k=0}u_{k}$ where $u_{k}=f_{k}(|x|)\varphi_{k}(x)$, we calculate
\begin{eqnarray*}
\frac{1}{n\omega_{n}}\int_{R^n} \frac{|\Delta %%@
u_{k}|^2}{|x|^{2m}}dx&=&\int_{0}^{R}r^{n-2m-1}(f''_{k}(r))^2dr+[(n-1)(2m+1)+2c_{k}]\int^{R}_{0}r^{%%@
n-2m-3}(f_{k}')^2dr\\
&+&c_{k}[c_{k}+(n-2m-4)(2m+2)]\int_{0}^{R}r^{n-2m-5}(f_{k}(r))^2dr\\
&\geq&(\frac{n+2m}{2})^2\int^{R}_{0}r^{n-2m-3}(f_{k}')^2dr+ %%@
\beta(|x|^{2m}; R)\int^{R}_{0}r^{n-1}(f_{k}')^2dr\\
&+&c_{k}\int^{R}_{0}r^{n-2m-3}(f_{k}')^2dr\\
&\geq& \frac{(n+2m)^{2}(n-2m-4)^2}{16}\int^{R}_{0}r^{n-2m-5}(f_{k})^2dr\\
&&+\beta 
(W; R)\frac{(n+2m)^2}{4}
\int^{R}_{0}W(r)r^{n-2m-3}(f_{k})^2dr\\
&+&\beta(|x|^{2m}; R)\int^{R}_{0}r^{n-1}(f_{k}')^2dr+c_{k}\beta(|x|^{2m}; R)\int^{R}_{0}r^{n-3}(f_{k})^2
dr\\
&=&\frac{(n+2m)^{2}(n-2m-4)^2}{16n\omega_{n}}\int_{R^n}\frac{u_{k}^2}{|x|^{2m+4}}dx\\
&+&\frac{\beta (W; R)}{n\omega_{n}}(\frac{(n+2m)^2}{4})
\int_{R^n}\frac{W(x)}{|x|^{2m+2}}u_{k}^2 dx+ \beta(|x|^{2m}; R)||u_{k}||_{W_{0}^{1,2}}.
\end{eqnarray*}
Hence (\ref{hrs}) holds. \hfill $\Box$

 We note that even for $m=0$ and $n\geq 4$, Theorem \ref{hrs-in} improves considerably Theorem %%@
A.2. in \cite{AGS}.

\section{Higher order Rellich inequalities}
In this section we will repeat the results obtained in the previous section to derive higher order %%@
Rellich inequalities with corresponding improvements. Let $W$ be a Bessel potential, $\beta_{n,m}$ %%@
be defined as in Theorem \ref{gen.hr} and 
\[\sigma_{n,m}=\beta (W; R)(\frac{(n+2m)^2}{4}+\frac{(n-2
m-\lambda-2)^2}{4}).\]
For the sake of convenience we make the following convention: $\prod\limits_{i=1}^{0}a_{i}=1.$
\begin{theorem}\label{h.o.rellich1}
Let $B_{R}$ be a ball of radius $R$ and $W$ be a Bessel potential on $B_{R}$ such that %%@
$\frac{W(r)}{W_r(r)}=-\frac{\lambda}{r}+f(r)$, where $f(r)\geq 0$ and $\lim_{r \rightarrow %%@
0}rf(r)=0$. Assume $m \in N$, $1\leq l\leq m$, and $2k+4m\leq n$. Then the following inequality %%@
holds for all $u \in C^{\infty}_{0}(B_R)$
\begin{equation}
\int_{B_{R}}\frac{|\Delta^{m}u|^2}{|x|^{2k}}dx \geq %%@
\prod\limits_{i=0}^{l-1}\beta_{n,k+2i}\int_{B_{R}}\frac{|\Delta^{m-l}u|^2}{|x|^{2k+4l}}dx+
\sum\limits_{i=0}^{l-1
}\sigma_{n,k+2i}\prod\limits_{j=1}^{l-1}\beta_{n,k+2j-2}\int_{B_{R}}\frac{W(x)|\Delta^{m-i-1}u|^2}%%@
{|x
|^{2k+4i+2}}dx
\end{equation}
\end{theorem}
{\bf Proof:} Follows directly from theorem \ref{gen.hr}. \hfill $\Box$

\begin{theorem}\label{h.o.rellich2}
Let $B_{R}$ be a ball of radius $R$ and $W$ be a Bessel potential on $B_{R}$ such that %%@
$\frac{W(r)}{W_r(r)}=-\frac{\lambda}{r}+f(r)$, where $f(r)\geq 0$ and $\lim_{r \rightarrow %%@
0}rf(r)=0$. Assume $m \in N$, $1\leq l\leq m$, and $2k+4m+2\leq n$. Then the following inequality %%@
holds for all $u \in C^{\infty}_{0}(B_R)$
\begin{eqnarray}
\int_{B_{R}}\frac{|\nabla \Delta^{m}u|^2}{|x|^{2k}}dx &\geq& (\frac{n-2k-2}{2})^2 %%@
\prod\limits_{i=0}^{l-1}\beta_{n,k+2i+1}\int_{B_{R}}\frac{|\Delta^{m-l}u|^2}{|x|^{2k+4l+2}}dx %%@
\nonumber\\
&+&(\frac{n-2k-2}{2})^2 %%@
\sum\limits_{i=0}^{l-1}\sigma_{n,k+2i+1}\prod\limits_{j=1}^{l-1}\beta_{n,k+2j-1}\int_{B_{R}}
\frac{W
(x)
|\Delta^{m-i-1}u|^2}{|x|^{2k+4i+4}}dx\nonumber \\
&+&\beta(W; R)\int_{B_{R}}W(x)\frac{|\Delta^{m}u|^2}{|x|^{2k}}dx
\end{eqnarray}
\end{theorem}
{\bf Proof:} Follows directly from Theorem \ref{main-CKN} and the previous theorem. \hfill $\Box$
\begin{remark} For $k=0$ Theorems \ref{h.o.rellich1} and \ref{h.o.rellich2} include Theorem 1.9 in %%@
\cite{TZ} as a special case.
\end{remark}

\begin{theorem}\label{h.o.rellich3}
Let $B_{R}$ be a ball of radius $R$ and $W$ be a Bessel potential on $B_{R}$ such that %%@
$\frac{W(r)}{W_r(r)}=-\frac{\lambda}{r}+f(r)$, where $f(r)\geq 0$ and $\lim_{r \rightarrow %%@
0}rf(r)=0$. Assume $m \in N$, $1\leq l\leq m-1$, and $2k+4m\leq n$. Then the following inequality %%@
holds for all $u \in C^{\infty}_{0}(B_R)$
\begin{eqnarray}
\int_{B_{R}}\frac{|\Delta^{m}u|^2}{|x|^{2k}}dx &\geq& a_{n,k}(\frac{n-2k-4}{2})^2 %%@
\prod\limits_{i=0}^{l-1}\beta_{n,k+2i+2}\int_{B_{R}}\frac{|\Delta^{m-l-1}u|^2}{|x|^{2k+4l+4}}dx
\nonumber\\
&+&a_{n,k}(\frac{n-2k-4}{2})^2 %%@
\sum\limits_{i=0}^{l-1}\sigma_{n,k+2i+2}\prod\limits_{j=1}^{l-1}\beta_{n,k+2j}\int_{B_{R}}\frac{W(%%@
x)|\Delta^{m-i-2}u|^2}{|x|^{2k+4i+6}}dx\nonumber \\
&+&\beta(W; R)a_{n,k}\int_{B_{R}}W(x)\frac{|\Delta^{m-1}u|^2}{|x|^{2k+2}}dx+\beta(W; %%@
R)\int_{B_{R}}W(x
)
\frac{|\nabla \Delta^{m-1}u|^2}{|x|^{2k}}dx
\end{eqnarray}
where $a_{n,m}$ is defined in Theorem \ref{gm.hr}.
\end{theorem}
{\bf Proof:} Follows directly from Theorem \ref{gm.hr} and the previous theorem. \hfill $\Box$\\

The following improves Theorem 1.10 in \cite{TZ} in many ways, since it is assumed there that %%@
$l\leq \frac{-n+8+2\sqrt{n^2-n+1}}{12}$ and $4m<n$. Even for $k=0$, Theorem \ref{h.o.rellich4} %%@
below shows that we can drop the first condition and replace the second one by $4m\leq n$.

\begin{theorem}\label {h.o.rellich4}
Let $B_{R}$ be a ball of radius $R$ and $W$ be a Bessel potential on $B_{R}$ such that . Assume $m %%@
\in N$, $1\leq l\leq m$, and $2k+4m\leq n$. Then the following inequality holds for all $u \in %%@
C^{\infty}_{0}(B_R)$
\begin{eqnarray}\label{nuts}
\int_{B_{R}}\frac{|\Delta^{m}u|^2}{|x|^{2k}}dx &\geq& %%@
\prod\limits_{i=1}^{l}\frac{a_{_{n,k+2i-2}}(n-2k-4i)^2}{4}\int_{B_{R}}\frac{|\Delta^{m-l}u|^2}{|x|%%@
^{2
k+4l}}dx\\
&+&\beta(W;R)\sum\limits_{i=1}^{l}\prod\limits_{j=1}^{l-1}\frac{a_{_{n,k+2j-2}}(n-2k-4j)^2}{4}\int%%@
%%@
_{B_{R
}}W(x)\frac{|\nabla \Delta^{m-i}u|^2}{|x|^{2k+4i-4}}dx\nonumber\\
&+&\beta(W;R)\sum\limits_{i=1}^{l}a_{_{n,k+2i-2}}\prod\limits_{j=1}^{l-1}\frac{a_{_{n,k+2j-2}}(n-2%%@
%%@
k-4j)^2}{4
}\int_{B_{R}}W(x)\frac{|\Delta^{m-i}u|^2}{|x|^{2k+4i-2}}dx,\nonumber
\end{eqnarray}
where $a_{n,m}$ are the best constants in inequality (\ref{gm-hr}). 
\end{theorem}
{\bf Proof:} Follows directly from Theorem \ref{gm.hr}. \hfill $\Box$

\section{Appendix (A): The class of Bessel potentials}

The Bessel equation associated to a potential $W$
\begin{equation*}\label{Bessel}
\hbox{$(B_W)$ \hskip 150pt $y''+\frac{1}{r}y'+W(r)y=0$\hskip 150pt} 
\end{equation*}
 is central to all results revolving around the inequalities of Hardy and Hardy-Rellich type. We %%@
summarize in 
this appendix the various properties of these equations that were used throughout this paper. 
\begin{definition} We say that a non-negative real valued $C^1$-function is a {\it Bessel %%@
potential on $(0, R)$} if there exists $c>0$ such that the equation $(B_{cW})$ has a positive %%@
solution on $(0, R)$.

The class of Bessel potentials on $(0, R)$ will be denoted by  ${\cal B}(0, R)$. 
  \end{definition}
Note that  the change of variable  $z(s)=y(e^{-s})$ maps the equation $y''+\frac{1}{r}y'+W(r)y=0$ %%@
into  
\begin{equation}
\hbox{$(B'_W)$ \hskip 150pt $z''+e^{-2s}W(e^{-s})z(s)=0.$\hskip 150pt} 
\end{equation}

On the other hand, the change of variables $\psi  
(t)=\frac{-e^{-t}y'(e^{-t})}{y(e^{-t})}$ maps  it into  the nonlinear equation
\begin{equation}
\hbox{$(B''_W)$ \hskip 150pt $\psi'(t)+\psi^{2}(t)+e^{-2t}W(e^{-t})=0$. \hskip 150pt} 
\end{equation}
This will allow us to relate the existence of positive solutions of $(B_W)$ to the non-oscillatory %%@
behaviour of equations $(B'_W)$ and $(B''_W)$. 

The theory of sub/supersolutions --applied to $(B''_W)$ (See Wintner \cite{win1, win2, Har})-- %%@
already yields, that if $(B_W)$ has a positive solution on an interval $(0, R)$ for some %%@
non-negative potential $W\geq 0$, then   for any $W$ such that $0\leq V \leq W$, the equation %%@
$(B_V)$ has also a positive solution on $(0, R)$. 
This  leads to the definition of the {\it weight} of a potential $W\in  {\cal B}(0, R)$ as:
 \begin{equation}
\hbox{$\beta (W;  R)=\sup\{c>0; \, (B_{cW})$ has a positive solution on $(0, R)$\}.}
 \end{equation}
 The following is now straightforward.
  \begin{proposition}  1)\,  The class ${\cal B}(0, R)$ is a closed convex and solid subset of %%@
$C^1(0, R)$. 
 
 2)\,  For every $W\in {\cal B}(0, R)$, the equation 
 \begin{equation*}\label{Bessel}
\hbox{$(B_W)$ \hskip 150pt $y''+\frac{1}{r}y'+\beta(W; R)W(r)y=0$\hskip 150pt} 
\end{equation*}
has a positive solution on $(0, R)$. 
\end{proposition}

The following gives an integral criteria for Bessel potentials. 
 
 \begin{proposition}\label{integral}  
  Let $W$ be a positive  locally integrable function on $\R$.   
  \begin{enumerate}
\item If 
$\liminf\limits_{r\rightarrow 0} \ln(r)\int^{r}_{0} sW(s)ds>-\infty$, 
 then for every $R>0$, there exists $\alpha:=\alpha(R)>0$ such that the scaled function %%@
$W_\alpha(x):=\alpha^2W(\alpha x)$ is a Bessel potential on $(0, R)$. 
\item  If
$\hbox{$\lim\limits_{r\rightarrow 0} \ln(r)\int^{r}_{0} sW(s)ds=-\infty$, } $
then there are no $\alpha, c>0$,  for which  
$W_{\alpha,c}=cW(\alpha |x|)$ is a Bessel potential on $(0, R)$.

\end{enumerate}
 \end{proposition}
{\bf Proof:} This relies on well known results concerning  the existence of  
non-oscillatory solutions (i.e., those $z(s)$ such that $z(s)>0$ for $s>0$ sufficiently large) for  
the second order linear differential equations
\begin{equation}\label{OSI_ODE}
z''(s)+a(s)z(s)=0,
\end{equation}
where $a$ is a locally integrable function on $\R$. For these equations, the following integral  
criteria are available. We refer to  \cite{Har, Hua, win1,win2, won}) among others for proofs and %%@
related results.
  \begin{description}
\item i)\,   If $\hbox{$\limsup_{t\rightarrow \infty } t\int^{\infty}_{t}a(s)ds<\frac{1}{4}$, } $ %%@
then Eq. (\ref{OSI_ODE}) is non-oscillatory. 
\item ii)\,  If 
$\hbox{$\liminf_{t\rightarrow \infty } t\int^{\infty}_{t}a(s)ds>\frac{1}{4}$,}$
 then Eq. (\ref{OSI_ODE}) is oscillatory.
\end{description}
 It follows that if $\liminf\limits_{r\rightarrow 0} \ln(r)\int^{r}_{0} sW(s)ds>-\infty$ holds,  %%@
then there exists $\delta>0$ such that $(B_W)$ has a positive solution on $(0, \delta)$. An easy %%@
scaling argument then shows that there exists $\alpha>0$ such that $W_\alpha(x):=\alpha^2W(\alpha %%@
x)$ is a Bessel potential on $(0, R)$. The rest of the proof is similar.\hfill $\square$\\

We now exhibit a few explicit Bessel potentials and compute their weights.  We use the following %%@
notation.
 \begin{equation}\label{log-def}
\hbox{$log^{(1)}(.)=log(.)$ \ \ and \ \  $log^{(k)}(.)=log(log^{(k-1)}(.))$ \ \ for\ \ $k\geq 2$.}
\end{equation} 
and 
\begin{equation}\label{X-def}
X_{1}(t)=(1-\log(t))^{-1}, \quad  X_{k}(t)=X_{1}(X_{k-1}(t)) \ \ \ \ k=2,3, ... ,
\end{equation}

\begin{theorem} \label{Bessel.theorem}{\rm \bf Explicit  Bessel potentials}
 
 \begin{enumerate}
 \item $  W \equiv 0$ is a Bessel potential on $(0,  R)$ for any $R>0$. 
\item The  Bessel function $J_{0}$ is a positive solution for equation $(B_W)$ with $W\equiv 1$, %%@
on $(0, z_0)$, where $z_{0}=2.4048...$ is the first zero of $J_0$. Moreover, $z_{0}$ is larger %%@
than the first root of any other solution for  $(B_1)$. In other words,

 for every $R>0$, 
\begin{equation}
\beta (1;  R)=\frac{z_0^2}{R^2}.
\end{equation}
\item   If $a<2$, then there exists $R_a>0$ such that $W (r)=r^{-a}$ is  a Bessel potential on $( %%@
0, R_a)$. 
\item For each $k\geq 1$ and $\rho>R( e^{e^{e^{.^{.^{e(k-times)}}}}} )$, the equation %%@
$(B_{\frac{1}{4}W_{k, \rho}})$ corresponding to the potential
\[
\hbox{$W_{k, \rho}(r)=\Sigma_{j=1}^kU_{j}$ where 
$
U_j(r)=\frac{1}{r^{2}}\big(\prod^{j}_{i=1}log^{(i)}\frac{\rho}{r}\big)^{-2}$}
\]
has a positive solution on $(0, R)$ that is explicitly given by 
 $
 \varphi_{k, \rho}(r)=( \prod^{k}_{i=1}log^{(i)}\frac{\rho}{r})^{\frac{1}{2}}.$
  On the other hand, the equation $(B_{_{\frac{1}{4}W_{k, \rho}+\lambda U_k}})$ corresponding to %%@
the potential $\frac{1}{4}W_{k, \rho}+\lambda U_k$ has no positive solution for any $\lambda>0$. 
   In other words, $W_{k, \rho}$ is a Bessel potential on $(0, R)$ with 
  \begin{equation}
\hbox{$ \beta (W_{k;  \rho}, R)=\frac{1}{4}$ for any $k\geq 1$.} 
 \end{equation}
 \item For each $k\geq 1$ and $R>0$, the equation $(B_{\frac{1}{4}\tilde W_{k, R}})$ corresponding %%@
to the potential 
\[
\hbox{$ \tilde W_{k, R}(r)=\Sigma_{j=1}^k\tilde U_{j}$ where  $
  \tilde U_j(r)=\frac{1}{r^{2}}X^{2}_{1}(\frac{r}{R})X^{2}_{2}(\frac{r}{R}) \ldots %%@
X^{2}_{j-1}(\frac{r}{R})X^{2}_{j}(\frac{r}{R})$}
\]
 has a positive solution on $(0, R)$ that is explicitly given by 
\begin{eqnarray*}
\varphi_{k}(r)=(X_{1}(\frac{r}{R})X_{2}(\frac{r}{R}) \ldots %%@
X_{k-1}(\frac{r}{R})X_{k}(\frac{r}{R}))^{-\frac{1}{2}}.
\end{eqnarray*}
On the other hand, the equation $(B_{\frac{1}{4}\tilde W_{k, R}+\lambda \tilde U_k})$ %%@
corresponding to the potential $\frac{1}{4}\tilde W_{k, R}+\lambda \tilde U_k$ has no positive %%@
solution for any $\lambda>0$.    In other words, $ \tilde W_{k, R}$ is a Bessel potential on $(0, %%@
R)$ with 
  \begin{equation}
\hbox{$\beta ({\tilde W}_{k, R};  R)=\frac{1}{4}$\quad  for any $k\geq 1$.} 
 \end{equation}

\end{enumerate}

\end{theorem} 
 
{\bf Proof:} 1) It is clear that $\phi (r)=-log(\frac{e}{R}r)$ is a positive solution of $(B_0)$ %%@
on $(0, R)$ for any $R>0$. 

2)The best constant for which the equation $y''+\frac{1}{r}y'+cy=0$ has a positive solution on %%@
$(0,R)$ is $\frac{z^{2}_{0}}{R^2}$, where $z_{0}=2.4048...$ is the first zero of Bessel function %%@
$J_{0}(z)$. Indeed if $\alpha$ is the first root of the an arbitrary solution of the Bessel %%@
equation $y''+\frac{y'}{r}+y(r)=0$, then we have $\alpha \leq z_{0}$. To see this let %%@
$x(t)=aJ_{0}(t)+bY_{0}(t)$, where $J_{0}$ and $Y_{0}$ are the two standard linearly independent %%@
solutions of Bessel equation, and $a$ and  $b$ are constants. Assume the first zero of $x(t)$ is %%@
larger than $z_{0}$.  Since the first zero of $Y_{0}$ is smaller than $z_{0}$, we have $a\geq0$. %%@
Also $b\leq0$, because $Y_{0}(t)\rightarrow -\infty$ as $t \rightarrow 0$. Finally note that %%@
$Y_{0}(z_{0})>0$, so if $b<0$, then $x(z_{0}+\epsilon)<0$  for $\epsilon$ sufficiently small. %%@
Therefore,  $b=0$ which is a contradiction.   

3) follows directly from the integral criteria. 

4) That $\phi_k$ is an explicit solution of the equation $(B_{\frac{1}{4}W_k})$ is %%@
straightforward. 
Assume now that there exists a positive function $\varphi$ such that
\begin{equation*}
-\frac{\varphi'(r)+r\varphi''(r)}{\varphi(r)}=\frac{1}{4}\sum^{k-1}_{j=1}\frac{1}{r}\big( %%@
\prod^{j}_{i=1}log^{(i)}\frac{\rho}{r}\big)^{-2}+(\frac{1}{4}+\lambda)\frac{1}{r}\big( %%@
\prod^{k}_{i=1}log^{(i)}\frac{\rho}{r}\big)^{-2}.
\end{equation*}
Define $f(r)=\frac{\varphi(r)}{\varphi_{k}(r)}>0$, and calculate, 
\[\frac{\varphi'(r)+r\varphi''(r)}{\varphi(r)}=\frac{\varphi_{k}'(r)+r\varphi_{k}''(r)}{\varphi_{k
}(r)}+\frac{f'(r)+rf''(r)}{f(r)}-\frac{f'(r)}{f(r)}\sum_{i=1}^{k}\frac{1}{\prod^{i}_{j=1}\log^{j}(%%@
\frac{\rho}{r})}.\]
Thus,
\begin{equation}\label{main-rel}
\frac{f'(r)+rf''(r)}{f(r)}-\frac{f'(r)}{f(r)}\sum_{i=1}^{k}\frac{1}{\prod^{i}_{j=1}\log^{j}(\frac{%%@
\rho}{r})}=-\lambda\frac{1}{r}\big( \prod^{k}_{i=1}log^{(i)}\frac{\rho}{r}\big)^{-2}.
\end{equation}
If now $f'(\alpha_{n})=0$ for some sequence $\{\alpha_{n}\}^{\infty}_{n=1}$ that converges to %%@
zero, then there exists a sequence $\{\beta_{n}\}^{\infty}_{n=1}$ that also converges to zero, %%@
such that  $f''(\beta_{n})=0$, and $f'(\beta_{n})>0$. But this contradicts (\ref{main-rel}), which %%@
means that $f$ is eventually monotone for $r$ small enough. We consider the two cases according to %%@
whether $f$ is increasing or decreasing:\\ \\
Case I:  Assume $f'(r)>0$ for $r>0$ sufficiently small. Then we will have
\[\frac{(rf'(r))'}{rf'(r)}\leq \sum_{i=1}^{k}\frac{1}{r\prod^{i}_{j=1}\log^{j}(\frac{\rho}{r})}.\]
Integrating once we get
\[f'(r)\geq \frac{c}{r\prod^{k}_{j=1}\log^{j}(\frac{\rho}{r})},\]
for some $c>0$.
Hence, $\lim_{r\rightarrow 0}f(r)=-\infty$ which is a contradiction. \\ \\
Case II: Assume $f'(r)<0$ for $r>0$ sufficiently small. Then
\[\frac{(rf'(r))'}{rf'(r)}\geq \sum_{i=1}^{k}\frac{1}{r\prod^{i}_{j=1}\log^{j}(\frac{\rho}{r})}. %%@
\]
Thus,
\begin{equation}\label{estim}
f'(r)\geq- \frac{c}{r\prod^{k}_{j=1}\log^{j}(\frac{\rho}{r})},
\end{equation}
for some $c>0$ and $r>0$ sufficiently small. On the other hand
\[\frac{f'(r)+rf''(r)}{f(r)}\leq-\lambda \sum^{k}_{j=1} \frac{1}{r}\big( %%@
\prod^{j}_{i=1}log^{(i)}\frac{R}{r}\big)^{-2}\leq %%@
-\lambda(\frac{1}{\prod^{k}_{j=1}\log^{j}(\frac{\rho}{r})})'.\]
Since $f'(r)<0$, there exists $l$ such that $f(r)>l>0$ for $r>0$ sufficiently small. From the %%@
above inequality we then have
\[bf'(b)-af'(a)<-\lambda %%@
l(\frac{1}{\prod^{k}_{j=1}\log^{j}(\frac{\rho}{b})}-\frac{1}{\prod^{k}_{j=1}\log^{j}(\frac{\rho}{a
})}).\]
From (\ref{estim}) we have $\lim_{a \rightarrow 0}af'(a)=0$. Hence, 
\[bf'(b)<-\frac{\lambda l }{\prod^{k}_{j=1}\log^{j}(\frac{\rho}{b})},\]
for every $b>0$, and
\[f'(r)<- \frac{\lambda l }{r\prod^{k}_{j=1}\log^{j}(\frac{\rho}{r})},\]
for $r>0$ sufficiently small. Therefore, 
\[\lim_{r\rightarrow 0}f(r)=+\infty,\]
and by choosing $l$ large enouph (e.g.,  $l>\frac{c}{\lambda})$ we get to contradict  %%@
$(\ref{estim})$.

The proof of 5) is similar and is left to the interested reader.\hfill $\square$

\section{Appendix (B): The evaluation of $a_{n,m}$}
Here we evaluate the best constants $a_{n,m}$ which appear in Theorem \ref{gm.hr}. 
\begin{theorem} Suppose $n\geq 1$ and $m\leq \frac{n-2}{2}$. Then for any $R>0$,   the constants
\[
a_{n,m}=\inf \left\{\frac{\int_{B_{R}}\frac{|\Delta u|^2}{|x|^{2m}}dx}{\int_{B_{R}}\frac{|\nabla %%@
u|^2}{|x|^{2m+2}}dx};\, u\in C^{\infty}_{0}(B_{R})\setminus \{0\}\right\}
\]
are given by the following expressions.
\begin{enumerate}
\item For $n=1$
\begin{itemize}
\item if $m \in (-\infty,-\frac{3}{2})\cup[-\frac{7}{6},-\frac{1}{2}]$, then
\[a_{1,m}=(\frac{1+2m}{2})^2\]
\item if $-\frac{3}{2}<m<-\frac{7}{6}$, then 
\[a_{1,m}=\min \{ %%@
(\frac{n+2m}{2})^2,\frac{(\frac{(n-4-2m)(n+2m)}{4}+2)^2}{(\frac{n-4-2m}{2})^2+2}\}.\]
\end{itemize} 
\item If $m=\frac{n-4}{2}$, then
\[a_{m,n}=\min\{(n-2)^2,n-1\}.\]
\item If $n\geq 2$ and $m\leq \frac{-(n+4)+2\sqrt{n^2-n+1}}{6}$, 
then $a_{n,m}=(\frac{n+2m}{2})^2$.
\item If $2 \leq n\leq 3$ and $\frac{-(n+4)+2\sqrt{n^2-n+1}}{6}<m\leq \frac{n-2}{2}$, or $n\geq 4$ %%@
and $\frac{n-4}{2}<m\leq \frac{n-2}{2}$, then 
\[a_{n,m}=\frac{(\frac{(n-4-2m)(n+2m)}{4}+n-1)^2}{(\frac{n-4-2m}{2})^2+n-1}.\]
\item For $n\geq 4$ and $\frac{-(n+4)+2\sqrt{n^2-n+1}}{6} < m< \frac{n-4}{2}$, define
$k^{*}=[(\frac{\sqrt{3}}{3}-\frac{1}{2})(n-2)].$
\begin{trivlist}
\item $\bullet$ If $k^{*}\leq 1$, then
\[a_{n,m}=\frac{(\frac{(n-4-2m)(n+2m)}{4}+n-1)^2}{(\frac{n-4-2m}{2})^2+n-1}.\]
\item $\bullet$ For $k^{*}>1$ the interval $(m^{1}_{0}:=\frac{-(n+4)+2\sqrt{n^2-n+1}}{6},m^{2}_{0}:= 
\frac{n-4}{2})$ can be divided in $2k^{*}-1$ subintervals. For $1 \leq k \leq k^{*}$ define
\[m^{1}_{k}:=\frac{2(n-5)-\sqrt{(n-2)^2-12k(k+n-2)}}{6},\]
\[m^{2}_{k}:=\frac{2(n-5)+\sqrt{(n-2)^2-12k(k+n-2)}}{6}.\]
If $m \in (m^{1}_{0},m^{1}_{1}]\cup[m^{2}_{1},m^{2}_{0})]$, then
\[a_{n,m}=\frac{(\frac{(n-4-2m)(n+2m)}{4}+n-1)^2}{(\frac{n-4-2m}{2})^2+n-1}.\]

\item $\bullet$ For $k\geq 1$ and $m \in (m^{1}_{k},m^{1}_{k+1}]\cup[m^{2}_{k+1},m^{2}_{k})$, then
\[a_{n,m}=\min\{\frac{(\frac{(n-4-2m)(n+2m)}{4}+k(n+k-2))^2}{(\frac{n-4-2m}{2})^2+k(n+k-2)},\frac{%%@
(\frac{(n-4-2m)(n+2m
)}{4}+(k+1)(n+k-1))^2}{(\frac{n-4-2m}{2})^2+(k+1)(n+k-1)}\}.\]
For $m \in (m^{1}_{k^{*}},m^{2}_{k^{*}})$, then
\[a_{n,m}=\min\{\frac{(\frac{(n-4-2m)(n+2m)}{4}+k^{*}(n+k^{*}-2))^2}{(\frac{n-4-2m}{2})^2+k^{*}(n+%%@
k^{*}-2)},\frac{(\frac{(n-4-2m)(n+2m)}{4}+(k^{*}+1)(n+k^{*}-1))^2}{(\frac{n-4-2m}{2})^2+(k^{*}+1)(%%@
n+k^{*}-1)}\}.\]
\end{trivlist}
\end{enumerate}
\end{theorem}

{\bf Proof:}  Letting %%@
$V(r)=r^{-2m}$ then,
\[W(r)-\frac{2V(r)}{r^2}+\frac{2V_r(r)}{r}-V_{rr}(r)= %%@
((\frac{n-2m-2}{2})^2-2-4m-2m(2m+1))r^{-2m-2}.\]
In order to satisfy condition (\ref{main.con}) we should have
\begin{equation}\label{exmain.con}
\frac{-(n+4)+2\sqrt{n^2-n+1}}{6} \leq m\leq \frac{-(n+4)+2\sqrt{n^2-n+1}}{6}.
\end{equation}
So, by Theorem \ref{main.hr} under the above condition we have $a_{n,m}=(\frac{n+2m}{2})^2$ as in the radial case. \\
For %%@
the rest of the proof we will use an argument similar to that of Theorem 6.4 in \cite{TZ} who computed $a_{n, m}$ in the case where $n\geq 5$ and for certain intervals of $m$.\\ %%@
Decomposing again $u \in C^{\infty}_{0}(B_{R})$ into spherical harmonics; $u=\Sigma^{\infty}_{k=0}u_{k}$, %%@
where $u_{k}=f_{k}(|x|)\varphi_{k}(x)$, one has
\begin{eqnarray}\label{N1}
\int_{\R^n}\frac{|\Delta %%@
u_{k}|^2}{|x|^{2m}}dx&=&\int_{\R^n}|x|^{-2m}(f''_{k}(|x|))^2dx+\left((n-1)(2m+1)+2c_{k}\right) %%@
\int_{\R^n}|x|^{-2m-2}(f_{k}')^2dx\\
&+&c_{k}(c_{k}+(n-4-2m)(2m+2))\int_{\R^n}|x|^{-2m-4}(f_{k})^2dx, \nonumber
\end{eqnarray}
\begin{equation}\label{N2}
\int_{\R^n}\frac{|\nabla u_{k}|^2}{|x|^{2m+2}}dx=\int_{\R^n}|x|^{-2m-2}(f_{k}')^2 %%@
dx+c_{k}\int_{\R^n}|x|^{-2m-4}(f_{k})^2dx.
\end{equation}
One can then prove as in \cite{TZ} that 
\begin{equation}\label{a-nm-k}
a_{n,m}=\min \left\{A(k,m,n); \, k\in \N \right\} 
\end{equation}
where
\begin{equation}
\hbox{$A(k,m,n)=\frac{(\frac{(n-4-2m)(n+2m)}{4}+c_{k})^2}{(\frac{n-4-2m}{2})^2+c_{k}}$  if  $m = \frac{n-4}{2}$}
\end{equation}
and 
 \begin{equation}
\hbox{$A(k,m,n):= c_{k}$ if  $m = \frac{n-4}{2}$ and $n+k>2$.}
\end{equation}
Note that when  $m=\frac{n-4}{2}$ and $n+k>2$, then $c_k\neq 0$.  Actually, this also holds for $n+k\leq 2$, in which case one deduces that if 
%First we prove that best constants $a_{n,m}$ are less than what theorem asserts.
%It follows from an argument similar to that of Theorem 4.6 in \cite{TZ} that
%\begin{eqnarray*}
%a_{n,m}&\leq& A(k,m,n):=\frac{(\frac{(n-4-2m)(n+2m)}{4}+c_{k})^2}{(\frac{n-4-2m}{2})^2+c_{k}} \ \ %%@
%if \ \ m \neq \frac{n-4}{2}\\
%a_{n,m}&\leq& c_{k}, \ \ if \ \ m = \frac{n-4}{2} \ \ and \ \ n+k>2.
%\end{eqnarray*}
%Hence by (\ref{N1}) and (\ref{N2}) we see that if 
$m = \frac{n-4}{2}$, then 
\[a_{n,m}= \min  
\{(n-2)^2=(\frac{n+2m}{2})^2, (n-1)=c_{1}\}
\] which is statement  2).

The rest of the proof consists of computing the infimum especially in the cases not considered in \cite{TZ}. 
For that we consider the function
\[f(x)=\frac{(\frac{(n-4-2m)(n+2m)}{4}+x)^2}{(\frac{n-4-2m}{2})^2+x}.\]
It is easy to check that $f'(x)=0$ at $x_{1}$ and $x_{2}$, where
\begin{eqnarray}
x_{1}&=&-\frac{(n-4-2m)(n+2m)}{4}\\
x_{2}&=&\frac{(n-4-2m)(-n+6m+8)}{4}.
\end{eqnarray}
Observe that for for $n\geq 2$, $\frac{n-8}{6}\leq \frac{n-4}{2}$. Hence, for $m\leq %%@
\frac{n-8}{6}$ both $x_{1}$ and $x_{2}$ are negative and hence $a_{n,m}= (\frac{n+2m}{2})^2$. %%@
Also note that 
\[\frac{-(n+4)-2\sqrt{n^2-n+1}}{6}\leq \frac{n-8}{6}\ \ for \ \ all \ \ n\geq 1.\] 
Hence, under the condition in 3) we have $a_{n,m}= (\frac{n+2m}{2})^2.$\\
 Also for $n=1$ if %%@
$m\leq -\frac{3}{2}$ both critical points are negative and we have $a_{1,m}\leq %%@
(\frac{1+2m}{2})^2$. Comparing $A(0,m,n)$ and $A(1,m,n)$ we see that $A(1,m,n)\geq A(0,m,n)$ if %%@
and only if (\ref{exmain.con}) holds. 

For $n=1$ and $-\frac{3}{2}<m<-\frac{7}{6}$ both $x_{1}$ and $x_{2}$ are positive. Consider the %%@
equations
\[x(x-1)=x_{1}=\frac{(2m+3)(2m+1)}{4},\]
and
\[x(x-1)=x_{2}=-\frac{(2m+3)(6m+7)}{4}.\]
By simple calculations we can see that all four solutions of the above two equations are less that %%@
two. Since, $A(1,m,1)<A(0,m,1)$ for $m<-\frac{7}{6}$, we have $a_{1,m}\leq \min %%@
\{A(1,m,1),A(2,m,1)\}$ and $1)$ follows.

For $n\geq 2$ and $\frac{n-4}{2}<m<\frac{n-2}{2}$ we have $x_{1}>0$ and $x_{2}<0$. Consider the %%@
equation
\[x(x+n-2)=x_{1}=-\frac{(n-4-2m)(n+2m)}{4}.\]
Then $\frac{2m+4-n}{2}$ and $-\frac{(2m+n)}{2}$ are solutions of the above equation and both are %%@
less than one. Since, for $n\geq 4$ 
\[\frac{n-2}{2}>\frac{-(n+4)+2\sqrt{n^2-n+1}}{6},\]
and $A(1,m,n)\leq A(0,m,n)$ for $m\geq \frac{-(n+4)+2\sqrt{n^2-n+1}}{6} $, the best constant is %%@
equal to what 4) claims. \\
$5)$ follows from an argument similar to that of Theorem 6.4 in %%@
\cite{TZ}.

%To prove that $a_{n,m}$ is the best constant, let $k$ be such that 
%\begin{eqnarray}\label{a-nm-k}
%a_{n,m}&=&\frac{(\frac{(n-4-2m)(n+2m)}{4}+k(n+k-2))^2}{(\frac{n-4-2m}{2})^2+k(n+k-2)} \ \ for \ \ %%@
%m\neq \frac{n-4}{2}\\
%a_{n,m}&=c_{k}& \ \ for \ \ m=\frac{n-4}{2},
%\end{eqnarray}
%and $k=0$ if $a_{n,m}=(n-2)^2$ and $m=\frac{n-4}{2}$. Set 
%\[u=|x|^{-\frac{n-4}{2}+m+\epsilon}\varphi_{k}(x)\varphi(|x|),\]
%where $\varphi_{k}(x)$ is an eigenfunction corresponding to the eigenvalue $c_{k}$ and %%@
%$\varphi(r)$ is a smooth cutoff function, such that $0 \leq \varphi \leq 1$, with $\varphi\equiv %%@
%1$ in $[0,\frac{1}{2}]$. We have
%\[\frac{\int_{B_{R}}\frac{|\Delta u|^2}{|x|^{2m}}dx}{\int_{B_{R}}\frac{|\nabla %%@
%u|^2}{|x|^{2m+2}}dx}=\frac{(-\frac{(n+2m)(n-4-2m)}{4}-c_{k}+\epsilon(2+2m+\epsilon))^2+O(1)}{(-
%\frac{n-4-2m}{2}+\epsilon)^2+c_{k}+O(1)}.\]
%Let now $\epsilon \rightarrow 0$ to obtain the result. Thus the inequality
%\[\int_{B_{R}}\frac{|\Delta u|^2}{|x|^{2m}}\geq a_{n,m}\int_{B_{R}}\frac{|\nabla %%@
%u|^2}{|x|^{2m+2}}dx,\]
%holds for all $u \in C^{\infty}_{0}(B_{R})$.  \hfill $\Box$

 \end{document}